\documentclass[10pt]{amsart}
\usepackage{amssymb}
\usepackage{amscd}

\theoremstyle{definition}
\newtheorem{thm}{Theorem}[section]
\newtheorem{lem}[thm]{Lemma}
\newtheorem{prp}[thm]{Proposition}
\newtheorem{dfn}[thm]{Definition}

\newtheorem{cnj}[thm]{Conjecture}

\newtheorem{rmk}[thm]{Remark}
\newtheorem{ntn}[thm]{Notation}
\newtheorem{exa}[thm]{Example}
\newtheorem{pbm}[thm]{Problem}

\newcommand{\beq}{\begin{equation}}
\newcommand{\eeq}{\end{equation}}
\newcommand{\beqr}{\begin{eqnarray*}}
\newcommand{\eeqr}{\end{eqnarray*}}
\newcommand{\bal}{\begin{align*}}
\newcommand{\eal}{\end{align*}}
\newcommand{\bei}{\begin{itemize}}
\newcommand{\eei}{\end{itemize}}
\newcommand{\limi}[1]{\lim_{{#1} \to \infty}}

\newcommand{\af}{\alpha}
\newcommand{\bt}{\beta}
\newcommand{\gm}{\gamma}
\newcommand{\dt}{\delta}
\newcommand{\ep}{\varepsilon}

\newcommand{\et}{\eta}
\newcommand{\ch}{\chi}
\newcommand{\io}{\iota}

\newcommand{\ld}{\lambda}
\newcommand{\sm}{\sigma}

\newcommand{\ph}{\varphi}
\newcommand{\ps}{\psi}
\newcommand{\rh}{\rho}
\newcommand{\om}{\omega}
\newcommand{\ta}{\tau}

\newcommand{\Ld}{\Lambda}

\newcommand{\Q}{{\mathbb{Q}}}
\newcommand{\Z}{{\mathbb{Z}}}
\newcommand{\R}{{\mathbb{R}}}
\newcommand{\C}{{\mathbb{C}}}
\newcommand{\N}{{\mathbb{N}}}

\newcommand{\id}{{\mathrm{id}}}

\newcommand{\sa}{{\mathrm{sa}}}

\newcommand{\diag}{{\mathrm{diag}}}

\newcommand{\rank}{{\mathrm{rank}}}

\newcommand{\card}{{\mathrm{card}}}
\newcommand{\Aut}{{\mathrm{Aut}}}
\newcommand{\Ad}{{\mathrm{Ad}}}

\newcommand{\Zq}[1]{{\Z_{#1}}}
\newcommand{\Zqt}{\Zq{2}}

\newcommand{\Zqn}{\Zq{n}}
\newcommand{\Cs}[3]{C^* (\Zq{#1}, #2, #3)}
\newcommand{\Csw}[3]{C^* (\Zq{#1}, \, #2, \, #3)}
\newcommand{\CZnAa}{\Cs{n}{A}{\af}}

\newcommand{\dirlim}{\displaystyle \lim_{\longrightarrow}}
\newcommand{\Mi}{M_{\infty}}

\newcommand{\andeqn}{\,\,\,\,\,\, {\mbox{and}} \,\,\,\,\,\,}
\newcommand{\QED}{\rule{0.4em}{2ex}}
\newcommand{\ts}[1]{{\textstyle{#1}}}
\newcommand{\ds}[1]{{\displaystyle{#1}}}

\newcommand{\ssum}[1]{{\ts{ {\ds{\sum}}_{#1} }}}

\newcommand{\ca}{C*-algebra}

\newcommand{\suca}{simple unital C*-algebra}

\newcommand{\idsuca}{infinite dimensional simple unital C*-algebra}

\newcommand{\idssucaz}{separable infinite dimensional
  simple unital C*-algebra with tracial rank zero}

\newcommand{\ct}{continuous}
\newcommand{\pj}{projection}

\newcommand{\hm}{homomorphism}

\newcommand{\Wolog}{Without loss of generality}

\newcommand{\tfae}{the following are equivalent}
\newcommand{\ifo}{if and only if}
\newcommand{\mops}{mutually orthogonal \pj s}

\newcommand{\cfn}{continuous function}
\newcommand{\hsa}{hereditary subalgebra}

\newcommand{\mvnt}{Murray-von Neumann equivalent}

\newcommand{\tRp}{tracial Rokhlin property}
\newcommand{\sRp}{strict Rokhlin property}
\newcommand{\tar}{tracially approximately representable}
\newcommand{\sar}{strictly approximately representable}
\newcommand{\fd}{finite dimensional}

\renewcommand{\S}{\subset}

\newcommand{\SM}{\setminus}
\newcommand{\I}{\infty}

\title[Actions with the tracial Rokhlin property]{Finite
  cyclic group actions with the tracial Rokhlin property}

\author{N.~Christopher Phillips}

\date{21~Sept.\  2006}

\address{Department of Mathematics, University  of Oregon,
       Eugene OR 97403-1222, USA.}

\email[]{ncp@darkwing.uoregon.edu}

\thanks{Research partially supported by
   NSF grants DMS~0070776 and DMS~0302401.}

\begin{document}

\setcounter{section}{-1}

\begin{abstract}
We give examples of actions of $\Z / 2 \Z$ on AF~algebras
and AT~algebras
which demonstrate the differences between the (strict) Rokhlin property
and the \tRp,
and between (strict) approximate representability
and tracial approximate representability.
Specific results include the following.
We determine exactly when a product type action of $\Z / 2 \Z$ on
a UHF~algebra has the \tRp;
in particular, unlike for the \sRp,
every UHF~algebra admits such an action.
We prove that Blackadar's action of $\Z / 2 \Z$
on the $2^{\infty}$ UHF~algebra,
whose crossed product is not~AF because it has nontrivial $K_1$-group,
has the \tRp,
and we give an example of an action of $\Z / 2 \Z$
on a simple unital AF~algebra which has the \tRp\  %
and such that the $K_0$-group of the crossed product
has torsion.
In particular, the crossed product of a simple unital AF~algebra by
an action of $\Z / 2 \Z$ with the \tRp\  %
need not be~AF.
We give examples of a \tar\  action of $\Z / 2 \Z$
on a simple unital AF~algebra which is nontrivial on $K_0,$
and of a \tar\  action of $\Z / 2 \Z$
on a simple unital AT~algebra with real rank zero which
is nontrivial on $K_1.$
\end{abstract}

\maketitle

\section{Introduction}\label{Sec:Intro}

\indent
The \tRp\  for actions of finite groups on \ca s was
introduced in~\cite{PhtRp1a} for the purpose of proving
that every simple higher dimensional noncommutative torus
is an AT~algebra (done in~\cite{PhtRp2}),
and proving that certain crossed products of such algebras
by finite cyclic groups are AF~algebras (done in~\cite{ELP}).
The purpose of this paper is to provide other examples
of actions of finite cyclic groups with the \tRp\  %
on \ca s with tracial rank zero.
We demonstrate by example the differences between the
(strict) Rokhlin property and the \tRp,
and between (strict) approximate representability
(Definition~3.6(2) of~\cite{Iz})
and its tracial analog,
tracial approximate representability
(Definition~3.2 of~\cite{PhtRp1a}).
(To emphasize the distinction with their tracial analogs,
in this paper we refer to the \sRp\  and to
strict approximate representability.)
In~\cite{Bl0},
Blackadar constructed an action of $\Z / 2 \Z$ on the
$2^{\infty}$~UHF algebra
such that the crossed product has nontrivial $K_1$-group,
and is hence not~AF.
As one of our examples, we prove that this action has the \tRp.
Earlier, in one of the exercises (10.11.3)
of his book~\cite{Bl},
Blackadar gave an example of an order two automorphism of $K_0 (A)$
for a simple separable AF~algebra $A$ such that,
if this automorphism could be implemented by an order two
automorphism of $A,$
then the resulting crossed product by $\Z / 2 \Z$
would have torsion in $K_0$ or nontrivial $K_1.$
With a very slight modification of Blackadar's algebra
(we use $\Z \big[ \frac{1}{3} \big] \oplus \Z$ instead of
$\Z \big[ \frac{1}{2} \big] \oplus \Z$),
our examples include actions of $\Z / 2 \Z$ on this AF~algebra
with the \tRp\  %
such that $K_0$ of the crossed product has a summand
isomorphic to $\Z / 2^n \Z,$
and also actions with the \tRp\  %
such that $K_1$ of the crossed product is nonzero.

Our results give counterexamples to various strengthenings
of results in~\cite{PhtRp1a}.
In particular,
using the notation $\Zq{n}$ for $\Z / n \Z$:
\begin{itemize}
\item
Even on a UHF~algebra, an action of $\Zqn$ with the \tRp\  need not have
the \sRp, and in fact the crossed product by such an action
can have nontrivial $K_1$-group, so need not be~AF.
See Example~\ref{CAR4},
and, for a simple AF~algebra which is not~UHF, Example~\ref{E:NoTor}.
\item
The crossed product of a simple unital AF~algebra
by an action of $\Zqn$ with the \tRp\  %
can have torsion in its $K_0$-group---even if the
action is \tar.
See Example~\ref{E:Torsion}.
\item
If an action $\af$ of $\Zqn$ on a \suca\  $A$
with tracial rank zero
is \sar\  and has the \tRp,
it does not follow that the automorphisms of the dual action are
approximately inner---even if $A$ is UHF and $\af$ is
locally representable in the sense of Section~I of~\cite{HR2}.
See Example~\ref{CAR2}.
\item
If an action $\af$ of $\Zqn$ on a \suca\  $A$
with tracial rank zero
is \tar\  and has the \sRp,
it does not follow that the dual action has
the \sRp---even if $A$ is~AT.
Use the dual of the action in Example~\ref{CAR4}.
\item
A tracially approximately inner automorphism of a \suca\  $A$
with tracial rank zero
need not be trivial on $K_0 (A)$---even if $A$ is AF and
$\af$ generates an action of $\Zqn$  with the \sRp.
Use the dual of the action in Example~\ref{CAR2}.
\item
A tracially approximately inner automorphism of a \suca\  $A$
with tracial rank zero
need not be trivial on $K_1 (A)$---even if $A$ is AT and
$\af$ generates an action of $\Zqn$  with the \sRp.
Use the dual of the action in Example~\ref{CAR4}.
\item
There is an action $\af$ of $\Zqn$ on a simple AF algebra such that
such that $\CZnAa$ is again a simple AF algebra,
but such that $\af$ does not have the \tRp.
This can happen even when the dual action has the \sRp\  and
$\af$ is \sar,
even locally representable in the sense of Section~I of~\cite{HR2}.
See Example~\ref{CAR3}.
\item
The dual of an action with the \tRp\  need not have the \tRp,
even when both the original algebra and the crossed product are
simple AF algebras,
and even when the original action in fact
has the \sRp.
Use the dual of the action in Example~\ref{CAR3}.
\end{itemize}

Examples involving actions on simple \ca s which do not
have tracial rank zero,
many of them on simple \ca s with tracial rank one,
will appear in~\cite{OP3}.

The first section contains preliminaries.
In particular, we recall the most important definitions,
and we give a new criterion for an action to have the \tRp.
The second section contains three examples of product type
actions of $\Zqt$ on the $2^{\infty}$~UHF algebra.
Along the way, we give a criterion for exactly
when such a product type action has the \tRp.
The third section contains the proof that
Blackadar's action has the \tRp.
The fourth section gives the examples of the kind suggested
by 10.11.3 of~\cite{Bl}.
In the last section, we state some open problems.

This paper contains the results
in Section~12 of the unpublished preprint~\cite{PhW},
with some of the examples put in a more general framework,
as well as some of the material from Section~13.
Section~4 here is entirely new.

We use the notation $\Zq{n}$ for $\Z / n \Z$;
the $p$-adic integers will not appear in this paper.
If $A$ is a \ca\  and $\af \colon A \to A$ is an automorphism
such that $\af^n = \id_A,$ then we write
$\CZnAa$ for the crossed product of
$A$ by the action of $\Zqn$ generated by $\af.$
We write $A_{\sa}$ for the set of selfadjoint elements of a \ca\  $A.$
We write $p \precsim q$ to mean that the \pj\  $p$ is \mvnt\  to a
sub\pj\  of $q,$ and $p \sim q$ to mean that
$p$ is \mvnt\  to $q.$
Also, $[a, b]$ denotes the additive commutator $a b - b a.$
We denote by $T (A)$ the set of all tracial states on a unital \ca\  $A,$
equipped with the weak* topology.
For any element of $T (A),$
we use the same letter for its standard extension to $M_n (A)$
for arbitrary $n,$
and to $\Mi (A) = \bigcup_{n = 1}^{\infty} M_n (A)$ (no closure).

I am grateful to Masaki Izumi for discussions concerning the
(strict) Rokhlin property.
In particular, he suggested Example~\ref{CAR3}.
I would also like to thank Hiroyuki Osaka for carefully
reading earlier versions of this paper, and catching a number
of misprints and suggesting many improvements.

\section{The tracial Rokhlin property and tracial approximate
  representability}\label{Sec:Prelim}

We recall the properties we will use.
Since we will be almost exclusively concerned with
stably finite \ca s,
we give the simpler versions of the \tRp\  and of
tracial approximate representability that are valid for the
finite case.
The following is Definition~1.1 of~\cite{PhtRp1a}.
It is essentially Definition~3.1 of~\cite{Iz},
although that definition is stated in terms of central sequences.
See the discussion after Definition~1.1 of~\cite{PhtRp1a}
for the equivalence of the two definitions.
We furthermore call the condition the \sRp,
to emphasize the distinction with the \tRp.

\begin{dfn}\label{ERPDfn}
Let $A$ be a unital \ca,
and let $\af \colon G \to \Aut (A)$
be an action of a finite group $G$ on $A.$
We say that $\af$ has the
{\emph{strict Rokhlin property}} if for every finite set
$F \S A,$ and every $\ep > 0,$
there are \mops\  $e_g \in A$ for $g \in G$ such that:
\begin{enumerate}
\item\label{ERPDfn:1} %
$\| \af_g (e_h) - e_{g h} \| < \ep$ for all $g, h \in G.$
\item\label{ERPDfn:2} %
$\| e_g a - a e_g \| < \ep$ for all $g \in G$ and all $a \in F.$
\item\label{ERPDfn:3} %
$\sum_{g \in G} e_g = 1.$
\end{enumerate}
\end{dfn}

The following is a special case of Lemma~1.11 of~\cite{PhtRp1a}.
We specialize to cyclic groups because we only use cyclic groups
in this paper,
and to establish notation for the \pj s used in the \tRp\  for
finite cyclic groups.

\begin{lem}\label{TRPCond}
Let $A$ be a finite infinite dimensional simple unital \ca,
and let $\af \in \Aut (A)$ satisfy $\af^n = \id_A.$
The action of $\Zqn$ generated by $\af$ has the
\tRp\   \ifo\  for every finite set
$F \S A,$ every $\ep > 0,$
and every nonzero positive element $x \in A,$
there are \mops\  $e_0, e_1, \dots, e_{n - 1} \in A$ such that:
\newcounter{TmpEnumi}
\begin{enumerate}
\item\label{TRPCond:1} %
$\| \af (e_j) - e_{j + 1} \| < \ep$ for $0 \leq j \leq n - 1,$
where by convention we take the indices mod $n,$
that is, $e_n = e_0.$
\item\label{TRPCond:2} %
$\| e_j a - a e_j \| < \ep$ for $0 \leq j \leq n - 1$ and all $a \in F.$
\item\label{TRPCond:3} %
With $e = \sum_{j = 0}^{n - 1} e_j,$ the \pj\  $1 - e$ is \mvnt\  to a
\pj\  in the \hsa\  of $A$ generated by $x.$
\setcounter{TmpEnumi}{\value{enumi}}
\end{enumerate}
\end{lem}

The following definition
is equivalent to Definition~3.6(2) of~\cite{Iz}.
The definition of~\cite{Iz}
is actually stated in terms of central sequences,
but it is equivalent to our definition by
Lemma~3.1 of~\cite{PhtRp1a}.
Also, we use the term \sar\  to emphasize the distinction
with \tar.

\begin{dfn}\label{AppRep}
Let $A$ be a unital \ca.
An action $\af \colon G \to \Aut (A)$
of a finite abelian group $G$ on $A$
is {\emph{\sar}}
if for every finite set $F \S A$ and every $\ep > 0,$
there are unitaries $w_g \in A$ such that:
\begin{enumerate}
\item\label{AppRep:1} %
$\| \af_g (a) - w_g a w_g^* \| < \ep$
for all $a \in F$ and all $g \in G.$
\item\label{AppRep:2} %
$\| w_g w_h - w_{g h} \| < \ep$ for all $g, h \in G.$
\item\label{AppRep:3} %
$\| \af_g (w_h) - w_h \| < \ep$ for all $g, h \in G.$
\end{enumerate}
\end{dfn}

\begin{lem}[Lemma~3.3 of~\cite{PhtRp1a}]\label{TARForFinite}
Let $A$ be a finite \idsuca,
and let $\af \colon G \to \Aut (A)$
be an action of a finite abelian group $G$ on $A.$
Then $\af$ is \tar\  %
\ifo\  for every finite set $F \S A,$ every $\ep > 0,$
and every positive element $x \in A$ with $\| x \| = 1,$
there are a \pj\  $e \in A$
and unitaries $w_g \in e A e$ such that:
\begin{enumerate}
\item\label{TARForFinite:1} %
$\| e a - a e \| < \ep$ for all $a \in F.$
\item\label{TARForFinite:2} %
$\| \af_g (e a e) - w_g e a e w_g^* \| < \ep$
for all $a \in F$ and all $g \in G.$
\item\label{TARForFinite:3} %
$\| w_g w_h - w_{g h} \| < \ep$ for all $g, h \in G.$
\item\label{TARForFinite:4} %
$\| \af_g (w_h) - w_h \| < \ep$ for all $g, h \in G.$
\item\label{TARForFinite:5} %
$1 - e$ is \mvnt\  to a
\pj\  in the \hsa\  of $A$ generated by $x.$
\end{enumerate}
\end{lem}

We have the following duality between
the (tracial) Rokhlin property
and (tracial) approximate representability.

\begin{thm}[Lemma~3.8 of~\cite{Iz}]\label{RPDualToAppRep}
Let $A$ be a unital \ca,
and let $\af \colon G \to \Aut (A)$
be an action of a finite abelian group $G$ on $A.$
Then:
\begin{enumerate}
\item\label{RPDualToAppRep:1} %
$\af$ has the \sRp\  \ifo\  %
${\widehat{\af}}$ is \sar.
\item\label{RPDualToAppRep:2} %
$\af$ is \sar\  \ifo\  %
${\widehat{\af}}$ has the \sRp.
\end{enumerate}
\end{thm}

\begin{thm}[Theorem~3.12 of~\cite{PhtRp1a}]\label{TRPDualToTAppRep}
Let $A$ be an \idsuca,
and let $\af \colon G \to \Aut (A)$
be an action of a finite abelian group $G$ on $A$
such that $C^* (G, A, \af)$ is also simple.
Then:
\begin{enumerate}
\item\label{TRPDualToTAppRep:1} %
$\af$ has the \tRp\  \ifo\  %
${\widehat{\af}}$ is \tar.
\item\label{TRPDualToTAppRep:2} %
$\af$ is \tar\  \ifo\  %
${\widehat{\af}}$ has the \tRp.
\end{enumerate}
\end{thm}

\begin{rmk}\label{TARImpTAI}
We point out that if $A$ is an \idsuca,
and if $\af \colon G \to \Aut (A)$
is a \tar\  %
action of a finite abelian group $G$ on $A,$
then each $\af_g$ is tracially approximately inner
in the sense of Definition~5.1 of~\cite{PhtRp1a}.
\end{rmk}

We now give a useful criterion for an action to have the \tRp.

\begin{lem}\label{ARPFromPosElts}
Let $A$ be a \idssucaz.
Let $\af \colon G \to \Aut (A)$
be an action of a finite group $G$ on $A.$
Suppose that for every finite set $F \S A$ and every $\ep > 0$
there are positive elements
$a_g \in A$ for $g \in G,$ with $0 \leq a_g \leq 1,$ such that:
\begin{enumerate}
\item\label{ARPFromPosElts:1} %
$a_g a_h = 0$ for $g, h \in G$ with $g \neq h.$
\item\label{ARPFromPosElts:2} %
$\| \af_g (a_h) - a_{g h} \| < \ep$ for all $g, h \in G.$
\item\label{ARPFromPosElts:3} %
$\| a_g c - c a_g \| < \ep$ for all $g \in G$ and $c \in F.$
\item\label{ARPFromPosElts:4} %
$\ta \left( 1 - \sum_{g \in G} a_g \right) < \ep$
for every tracial state $\ta$ on $A.$
\end{enumerate}
Then $\af$ has the \tRp.
\end{lem}

\begin{proof}
We verify the criterion of Theorem~5.3 of~\cite{ELP}
using Conditions~(1), (2), and~(3$'$) there.
Thus,
let $S \S A$ be finite and let $\ep > 0.$
Set $n = \card (G).$
Choose $\ep_0 > 0$ so small that
$2 \sqrt{(2 n + 1) \ep_0} + (2 n + 5) \ep_0 \leq \ep.$

Apply the hypothesis
with $S$ as given and with $\ep_0$ in place of $\ep.$
For each $g \in G,$ use real rank zero
to choose a \pj\  $e_g \in {\overline{a_g A a_g}}$
such that
$\| e_g a_g - a_g \| < \ep_0$ and $\| a_g e_g - a_g \| < \ep_0.$
Then $\| e_g a_g e_g - a_g \| < 2 \ep_0.$
So
\[
\sum_{h \in G} \ta (e_h a_h e_h)
 > \sum_{h \in G} \ta (a_h) - 2 n \ep_0
 > 1 - (2 n + 1) \ep_0.
\]
Since $\sum_{h \in G} \ta (e_h) \leq 1$
(because the $e_h$ are orthogonal)
and $e_h a_h e_h \leq e_h$ for all $h \in G,$
we get, for each $g \in G$ and $\ta \in T (A),$
\[
\ta (e_g)
 \leq 1 - \sum_{h \neq g} \ta (e_h)
 \leq 1 - \sum_{h \neq g} \ta (e_h a_h e_h)
 < \ta (e_g a_g e_g) + (2 n + 1) \ep_0.
\]
Therefore $0 \leq \ta (e_g - e_g a_g e_g) < (2 n + 1) \ep_0.$

Recall the $L^2$-norm (or seminorm)
associated with a tracial state $\ta$ of a \ca\  $A,$ given by
$\| a \|_{2, \ta} = \ta (a^* a)^{1/2},$
as used in Section~5 of~\cite{ELP}.
Using $\| (e_g - e_g a_g e_g)^{1/2} \| \leq 1$ at the third step,
we now have
\begin{align*}
\| e_g - e_g a_g e_g \|_{2, \ta}
& \leq \big\| (e_g - e_g a_g e_g)^{1/2} \big\|_{2, \ta}
               \big\| (e_g - e_g a_g e_g)^{1/2} \big\|        \\
& \leq \ta (e_g - e_g a_g e_g)^{1/2}
  < \sqrt{(2 n + 1) \ep_0}.
\end{align*}
so that
\[
\| e_g - a_g \|_{2, \ta}
 \leq \| e_g - e_g a_g e_g \|_{2, \ta} + \| e_g a_g e_g - a_g \|
 < \sqrt{(2 n + 1) \ep_0} + 2 \ep_0.
\]

Now we can verify
Conditions~(1), (2), and~(3$'$) of Theorem~5.3 of~\cite{ELP}.
Let $\ta \in T (A).$
Then for $g, h \in G$ we have
\begin{align*}
\| \af_g (e_h) - e_{g h} \|_{2, \ta}
 & \leq \| e_g - a_g \|_{2, \ta} + \| e_g - a_g \|_{2, \ta}
       + \| \af_g (a_h) - a_{g h} \|          \\
 & < 2 \big( \sqrt{(2 n + 1) \ep_0} + 2 \ep_0 \big) + \ep_0
   \leq \ep,
\end{align*}
and for $g \in G$ and $c \in S$ we have
\[
\| [e_g, c] \|_{2, \ta}
  \leq 2 \| e_g - a_g \|_{2, \ta} + \| [a_g, c] \|
  < 2 \big( \sqrt{(2 n + 1) \ep_0} + 2 \ep_0 \big) + \ep_0
  \leq \ep.
\]
Also, using
\[
\ta (e_g)
  \geq \ta (e_g a_g e_g)
  \geq \ta (a_g) - \| a_g - e_g a_g e_g \|
  > \ta (a_g) - 2 \ep_0,
\]
we get
\[
\ta \left( 1 - \ssum{g \in G} e_g \right)
  < (2 n + 1) \ep_0
  \leq \ep.
\]
An application of Theorem~5.3 of~\cite{ELP}
completes the proof.
\end{proof}

\begin{thm}\label{ARPFromPosElts2}
Let $A$ be a \idssucaz.
Let $\af \colon G \to \Aut (A)$
be an action of a finite group $G$ on $A.$
Suppose that for every finite set $F \S A$ and every $\ep > 0$
there are positive elements
$a_g \in A$ for $g \in G,$ with $0 \leq a_g \leq 1,$ such that:
\begin{enumerate}
\item\label{ARPFromPosElts2:1} %
$\| a_g a_h \| < \ep$ for $g, h \in G$ with $g \neq h.$
\item\label{ARPFromPosElts2:2} %
$\| \af_g (a_h) - a_{g h} \| < \ep$ for all $g, h \in G.$
\item\label{ARPFromPosElts2:3} %
$\| a_g c - c a_g \| < \ep$ for all $g \in G$ and $c \in F.$
\item\label{ARPFromPosElts2:4} %
$\ta \left( 1 - \sum_{g \in G} a_g \right) < \ep$
for every tracial state $\ta$ on $A.$
\end{enumerate}
Then $\af$ has the \tRp.
\end{thm}

\begin{proof}
We prove that the hypotheses imply those of Lemma~\ref{ARPFromPosElts}.
Let $F \S A$ be finite and let $\ep > 0.$
\Wolog\  $\| c \| \leq 1$ for all $c \in F.$
Let $B$ be the universal (nonunital) \ca\  generated by
selfadjoint elements $d_g$ for $g \in G,$ subject to the
relations $0 \leq d_g \leq 1$ and $d_g d_h = 0$ for $g \neq h.$
Let $n = \card (G).$
Then $B$ is isomorphic to the cone over $\C^n,$ and is
hence a projective \ca.
(See Lemmas 8.1.3 and 10.1.5 and Theorem 10.1.11 of~\cite{Lr}.)
Therefore it is semiprojective (Definition 14.1.3 of~\cite{Lr}).
Thus, using Theorem 14.1.4 of~\cite{Lr}
(see Definition 14.1.1 of~\cite{Lr}, and take $B$ there to be $\{ 0 \}$),
there is $\dt > 0$ such that whenever $D$ is a \ca\  and
$d_g \in D,$ for $g \in G,$
are positive elements with $0 \leq d_g \leq 1$ and such that
$\| d_g d_h \| < \dt$ for $g \neq h,$ then there are
positive elements
$a_g \in D$ for $g \in G$ with $0 \leq a_g \leq 1$ such that
$a_g a_h = 0$ for $g \neq h$ and
\[
\| a_g - d_g \|
  < \min \left( \ts{ \frac{1}{4} } \ep, \,
                 \ts{ \frac{1}{2} } n^{-1} \ep \right)
\]
for all $g \in G.$
Apply our hypotheses with $F$ as given and with
$\min \left( \ts{ \frac{1}{2} } \ep, \, \dt \right)$ in place of $\ep,$
obtaining $d_g \in A$ for $g \in G,$
and let $a_g \in A$ be as above.
The relation $a_g a_h = 0$ for $g \neq h$
is Condition~(\ref{ARPFromPosElts:1}) of Lemma~\ref{ARPFromPosElts}.
For~(\ref{ARPFromPosElts:2}), estimate
\[
\| \af_g (a_h) - a_{g h} \|
  \leq \| \af_g (d_g) - d_{g h} \| + \| a_g - d_g \|
            + \| a_{g h} - d_{g h} \|
  < \ts{ \frac{1}{2} } \ep
            + \ts{ \frac{1}{4} } \ep + \ts{ \frac{1}{4} } \ep
  = \ep.
\]
For~(\ref{ARPFromPosElts:3}),
use $\| c \| \leq 1$ for $c \in F$ to estimate
\[
\| a_g c - c a_g \| \leq \| d_g c - c d_g \| + 2 \| a_g - d_g \|
     < \ts{ \frac{1}{2} } \ep + 2 \left( \ts{ \frac{1}{4} } \ep \right)
     = \ep.
\]
Condition~(\ref{ARPFromPosElts:4})
follows from $\| a_g - d_g \| < \frac{1}{2} n^{-1} \ep$
for $g \in G$
and $\ta \left( 1 - \sum_{g \in G} d_g \right) < \frac{1}{2} \ep.$
This completes the proof.
\end{proof}

\section{Product type actions}\label{Sec:Exs}

The three main examples in this section are product type actions
of $\Zqt$ on the $2^{\infty}$ UHF algebra $D.$
In each case we represent $D$ as an infinite tensor product,
and the automorphism generating the action as an infinite
tensor product of inner automorphisms.

It is useful to first develop some theory.
For simplicity, we consider only actions of $\Zqt.$
This class of actions has been considered several times before,
particularly in~\cite{FM} and~\cite{Ks0}
(which concentrate on the fixed point algebras)
and in~\cite{HR2}.
For the convenience of the reader,
and to establish notation, we start from scratch.
There is, however, some overlap with earlier work.
The formula in Lemma~\ref{Condense}(\ref{Condense:5}) is given,
in different terms, in~1.3 of~\cite{FM}.
Lemma~1.6.1 of~\cite{FM} gives the equivalence
of (\ref{CondForRP:1}) and~(\ref{CondForRP:3})
of Proposition~\ref{CondForRP},
together with several other
equivalent conditions similar to~(\ref{CondForRP:1}).
Example~III.4 of~\cite{HR2} contains the equivalence
of (\ref{CondForTRP:3}) and~(\ref{CondForTRP:4})
of Proposition~\ref{CondForTRP}, and the final statement.
The equivalence of (\ref{CondForOut:2}) and~(\ref{CondForOut:3})
in Proposition~\ref{CondForOut} is Lemma~1.4.1 of~\cite{FM}.

A complete description of the ordered K-theory of the crossed product
seems complicated, except when it is a UHF~algebra.
We say little about it here.
A partial description can be found in~2.2 of~\cite{FM}.
(It is actually intended for the fixed point algebra,
but the situation considered also covers crossed products.)
A different partial description can be found in Example~III.4
of~\cite{HR2}.

\begin{lem}\label{Z2Struct}
Let $D$ be an infinite tensor product \ca\  and let $\af \in \Aut (D)$
be an automorphism of order two, of the form
\[
D = \bigotimes_{n = 1}^{\I} M_{k (n)} \andeqn
\af = \bigotimes_{n = 1}^{\I}
 \Ad ( p_n - q_n ),
\]
with $k (n) \in \N$ and where $p_n, \, q_n \in M_{k (n)}$
are \pj s with $p_n + q_n = 1.$
Let $D_n = \bigotimes_{m = 1}^{n} M_{k (m)},$
so that $D = \dirlim D_n,$
and write $D_n = M_{t (n)},$ where $t (n) = \prod_{m = 1}^{n} k (m).$
Then the direct system of crossed
products can be identified as
\[
C^* (\Zqt, D_n) \cong M_{t (n)} \oplus M_{t (n)},
\]
with the dual action given by the flip $\sm_n (a, b) = (b, a)$
for $a, \, b \in M_{t (n)},$
and where the maps
\[
\ps_n \colon M_{t (n - 1)} \oplus M_{t (n - 1)}
         \to M_{t (n)} \oplus M_{t (n)}
\]
of the direct system
are given by
\[
\ps_n (a, b)
 = (a \otimes p_n + b \otimes q_n, \, b \otimes p_n + a \otimes q_n)
\]
for $a, \, b \in M_{t (n - 1)}.$
Moreover:
\begin{enumerate}
\item\label{Z2Struct:1} %
$\Cs{2}{D}{\af}$ is a unital AF~algebra.
\item\label{Z2Struct:2} %
$K_0 (\Cs{2}{D}{\af})$ is the direct limit of the system
\[
\Z^2 \stackrel{T_1}{\longrightarrow}
\Z^2 \stackrel{T_2}{\longrightarrow}
\Z^2 \stackrel{T_3}{\longrightarrow} \cdots,
\]
in which the $n$th term
is identified with $K_0 (M_{t (n)} \oplus M_{t (n)})$
in the obvious way,
with
\[
T_n = \left( \begin{array}{cc}  \rank (p_n)  & \rank (q_n)  \\
               \rank (q_n)  & \rank (p_n) \end{array} \right).
\]
The isomorphism identifies
${\widehat{\af}}_*$ with the direct limit of the
maps $(j, k) \mapsto (k, j)$ on $\Z^2.$
\item\label{Z2Struct:3} %
The action of $\Zqt$ generated by $\af$
is \sar\  (Definition~\ref{AppRep}).
\item\label{Z2Struct:4} %
The dual action on $\Cs{2}{D}{\af}$ has
the \sRp\   (Definition~\ref{ERPDfn}).
\item\label{Z2Struct:5} %
For any subset $S \subset \N,$
the automorphism $\af$ is unchanged if
$p_n$ and $q_n$ are exchanged for every $n \in S.$
\end{enumerate}
\end{lem}

\begin{proof}
For $n \in \N$ and a unitary $v \in M_n$ with $v^2 = 1,$
we use the isomorphism
$\Csw{2}{M_n}{\Ad (v)} \to M_n \oplus M_n$
which sends $a \in M_n$ to $(a, a)$ and the canonical unitary
of the crossed product to $(v, \, - v).$
The identification of the direct system of crossed products is then a
calculation, which we omit.
The K-theory computation in Part~(\ref{Z2Struct:2}) follows.
Parts~(\ref{Z2Struct:1}) and~(\ref{Z2Struct:3}) are immediate.
Part~(\ref{Z2Struct:4})
follows from Part~(\ref{Z2Struct:3})
and Theorem~\ref{RPDualToAppRep}(\ref{RPDualToAppRep:2}).
Part~(\ref{Z2Struct:5}) follows from the fact that $\Ad (-w) = \Ad (w).$
\end{proof}

\begin{ntn}\label{ProdTypeNtn}
For a product type automorphism $\af$ as in Lemma~\ref{Z2Struct},
we assume, as Lemma~\ref{Z2Struct}(\ref{Z2Struct:5}) says we may,
that $\rank (p_n) \geq \rank (q_n)$ for all $n.$
We then define
\[
\ld_n
  = \frac{\rank (p_n) - \rank (q_n)}{\rank (p_n) + \rank (q_n)}
  \geq 0
\]
for $n \in \N,$
and, for $m \leq n,$
\[
\Ld (m, n) = \ld_{m + 1} \ld_{m + 2} \cdots \ld_{n}
\andeqn
\Ld (m, \infty) = \limi{n} \Ld (m, n).
\]
We furthermore use $\af$ interchangeably for the action of $\Zqt$
and for the nontrivial automorphism which generates it,
and we use ${\widehat{\af}}$ interchangeably for the
dual action of $\Zqt$ on $\Cs{2}{D}{\af}$
and for the nontrivial automorphism which generates it.
\end{ntn}

It is useful to be able to condense the infinite tensor
product description by omitting intermediate terms.

\begin{lem}\label{Condense}
Let $D$ and $\af \in \Aut (D)$
be of the form
\[
D = \bigotimes_{n = 1}^{N} M_{k (n)} \andeqn
\af = \bigotimes_{n = 1}^{N} \Ad ( p_n - q_n ),
\]
with $k (n) \in \N$ and where $p_n, \, q_n \in M_{k (n)}$
are \pj s with $p_n + q_n = 1$
and $\rank (p_n) \geq \rank (q_n).$
Set $t = \prod_{n = 1}^{N} k (n).$
Then there exist \pj s $p, q \in M_t,$ with $p + q = 1$
and $\rank (p) \geq \rank (q),$
such that there is an isomorphism $D \cong M_t$
which intertwines $\af$ with $\Ad (p - q).$
Moreover: 
\begin{enumerate}
\item\label{Condense:1} %
With $T_n$ as in Lemma~\ref{Z2Struct}(\ref{Z2Struct:2}),
the ranks of $p$ and $q$ are determined
by the equation
\[
\left( \begin{array}{cc}  \rank (p)  & \rank (q)  \\
               \rank (q)  & \rank (p) \end{array} \right)
       = T_N T_{N - 1} \cdots T_1.
\]
\item\label{Condense:2} %
If for some $n$ we have $\rank (p_n) = \rank (q_n),$
then $\rank (p) = \rank (q).$
\item\label{Condense:3} %
If for all $n$ we have $\rank (p_n) > \rank (q_n),$
then $\rank (p) > \rank (q).$
\item\label{Condense:4} %
If for some $n$ we have $\rank (q_n) > 0,$
then $\rank (q) > 0.$
\item\label{Condense:5} %
We have
\[
\frac{\rank (p) - \rank (q)}{\rank (p) + \rank (q)}
  = \prod_{n = 1}^N
    \frac{\rank (p_n) - \rank (q_n)}{\rank (p_n) + \rank (q_n)}.
\]
\end{enumerate}
\end{lem}

\begin{proof}
For all parts,
the general case follows from the case $N = 2$ by induction,
so we do only that case.
Set
\[
p = p_1 \otimes p_2 + q_1 \otimes q_2
\andeqn
q = p_1 \otimes q_2 + q_1 \otimes p_2.
\]
Then one checks that $p + q = 1$ and
$(p_1 - q_1) \otimes (p_2 - q_2) = p - q.$
Taking the differences of the dimensions of the eigenspaces,
we get
\[
\rank (p) - \rank (q)
  = [\rank (p_1) - \rank (q_1)] [\rank (p_2) - \rank (q_2)].
\]
So $\rank (p) \geq \rank (q).$
We have now proved the main statement and Part~(\ref{Condense:1}).
Parts~(\ref{Condense:2}), (\ref{Condense:3}),
and~(\ref{Condense:4}) are immediate.
Part~(\ref{Condense:5}) follows by dividing the last equation by
$t = k (1) k (2).$
\end{proof}

The implication from~(\ref{CondForRP:3}) to~(\ref{CondForRP:1})
in the following lemma
is related to Example~3.2 of~\cite{Iz}.

\begin{prp}\label{CondForRP}
Let $\af \in \Aut (D)$ be a product type automorphism of order~$2$
as in Lemma~\ref{Z2Struct},
and adopt Notation~\ref{ProdTypeNtn}.
Then \tfae:
\begin{enumerate}
\item\label{CondForRP:1} %
The action $\af$ has the \sRp.
\item\label{CondForRP:2} %
The dual action ${\widehat{\af}}$ is \sar.
\item\label{CondForRP:3} %
There are infinitely many $n \in \N$
such that $\rank (p_n) = \rank (q_n)$
(that is, $\ld_n = 0$).
\item\label{CondForRP:4} %
$\Cs{2}{D}{\af}$ is a UHF~algebra.
\item\label{CondForRP:5} %
$K_0 (\Cs{2}{D}{\af})$ is totally ordered.
\item\label{CondForRP:6} %
${\widehat{\af}}_*$ is trivial on $K_0 (\Cs{2}{D}{\af}).$
\end{enumerate}
\end{prp}

\begin{proof}
That~(\ref{CondForRP:1}) implies~(\ref{CondForRP:2}) follows from
Theorem~\ref{RPDualToAppRep}(\ref{RPDualToAppRep:1}).
It is immediate that~(\ref{CondForRP:4}) implies~(\ref{CondForRP:5}).
Also, (\ref{CondForRP:2}) implies~(\ref{CondForRP:6})
because approximately inner automorphisms
are trivial on K-theory.

We prove that~(\ref{CondForRP:3})
implies (\ref{CondForRP:1}) and~(\ref{CondForRP:4}).
Using Lemma~\ref{Condense}(\ref{Condense:2}),
we may pass to a subsystem
and assume that $\rank (p_n) = \rank (q_n)$ for all $n.$
Then for every $n$ there is $v_n \in M_{k (n)}$ such that
$v_n v_n^* = p_n$ and $v_n^* v_n = q_n.$

To verify the \sRp, thus proving~(\ref{CondForRP:1}),
it suffices to consider finite subsets
$F$ of the dense subalgebra
$\bigcup_{n = 1}^{\I} D_n \S D.$
Then there is $N$ such that $F \S D_N.$
Define $f \in M_{k (N + 1)}$ by
$f = \frac{1}{2} (p_{N + 1} + q_{N + 1} + v_{N + 1} + v_{N + 1}^*).$
Computations show that $f$ is a \pj\  and
$(p_{N + 1} - q_{N + 1}) f (p_{N + 1} - q_{N + 1}) = 1 - f.$
The \pj s required for the \sRp\   using the given set
$F$ and any $\ep > 0$ are then
\[
e_0 = 1_{D_{N}} \otimes f,
\,\,  e_1 = 1_{D_{N}} \otimes (1 - f) \in D_{N + 1} \S D.
\]

To prove~(\ref{CondForRP:4}),
start by observing that $v_n + v_n^*$ is a unitary of order~$2$
such that $\Ad (v_n + v_n^*)$ exchanges $p_n$ and $q_n.$
So we can write the map $\ps_n$ of Lemma~\ref{Z2Struct}
as the composite of the maps
\[
\rh_n \colon M_{t (n - 1)} \oplus M_{t (n - 1)} \to M_{t (n)}
\andeqn
\sm_n \colon M_{t (n)} \to M_{t (n)} \oplus M_{t (n)}
\]
given by
$\rh_n (a, b) = a \otimes p_n + b \otimes q_n$
for $a, \, b \in M_{t (n - 1)}$
and
\[
\sm_n (a)
 = (a, \,\, [1 \otimes (v_n + v_n^*)] a [1 \otimes (v_n + v_n^*)])
\]
for $a \in M_{t (n)}.$
Then $D$ is the direct limit of the system
\[
M_{t (1)} \stackrel{\rh_2 \circ \sm_1}{\longrightarrow}
M_{t (2)} \stackrel{\rh_3 \circ \sm_2}{\longrightarrow}
M_{t (3)} \stackrel{\rh_4 \circ \sm_3}{\longrightarrow} \cdots,
\]
which is obviously a UHF~algebra.

We complete the proof by showing
that (\ref{CondForRP:5})~and~(\ref{CondForRP:6})
both imply~(\ref{CondForRP:3}),
which we do by contradiction.
Assume~(\ref{CondForRP:3}) fails.
Choose $N_0$ such that 
$\rank (p_n) > \rank (q_n)$ for all $n \geq N_0.$
Identifying $K_0 (C^* (\Zqt, M_{t (N_0)}))$ with $\Z^2$
as in Lemma~\ref{Z2Struct}(\ref{Z2Struct:2}),
set $\et_{N_0} = (1, \, -1) \in K_0 (C^* (\Zqt, M_{t (N_0)})).$
Let $N > N_0.$
Using Lemma~\ref{Condense}(\ref{Condense:1}),
we get
\[
T_N T_{N - 1} \cdots T_{N_0 + 1}
 = \left( \begin{array}{cc}  \rank (e_N)  & \rank (f_N)  \\
                \rank (f_N)  & \rank (e_N) \end{array} \right)
\]
for \pj s $e_N, f_N \in \bigotimes_{n = N_0 + 1}^{N} M_{k (n)}$
with $e_N + f_N = 1.$
Moreover,
$\rank (e_N) > \rank (f_N)$ by Lemma~\ref{Condense}(\ref{Condense:3}).
Lemma~\ref{Z2Struct}(\ref{Z2Struct:2})
implies that the image $\et_N$ of $\et_{N_0}$
in $K_0 (C^* (\Zqt, M_{t (N)}))$
is $(l_N, \, - l_N)$ with $l_N = \rank (e_N) - \rank (f_N).$
Since $l_N \neq 0$ for all $N,$
the image $\et$ of $\et_{N_0}$ in $K_0 (C^* (\Zqt, D, \af))$ is nonzero.
Lemma~\ref{Z2Struct}(\ref{Z2Struct:2})
implies that ${\widehat{\af}}_* (\et) = - \et.$
Since $K_0 (C^* (\Zqt, D, \af))$ has no torsion,
we have $- \et \neq \et,$
so~(\ref{CondForRP:6}) fails.
Moreover, for every $N \geq N_0,$ neither $\et_N > 0$ nor $- \et_N > 0.$
Therefore neither $\et > 0$ nor $- \et > 0,$
whence~(\ref{CondForRP:5}) fails.
\end{proof}

The last statement in the following proposition is a generalization
of a statement at the end of Example~3.14 of~\cite{Iz0}.
It is also in Example~III.4 of~\cite{HR2}.

\begin{prp}\label{CondForTRP}
Let $\af \in \Aut (D)$ be a product type automorphism of order~$2$
as in Lemma~\ref{Z2Struct},
and adopt Notation~\ref{ProdTypeNtn}.
Then \tfae:
\begin{enumerate}
\item\label{CondForTRP:1} %
The action $\af$ has the \tRp.
\item\label{CondForTRP:2} %
The dual action ${\widehat{\af}}$ is \tar.
\item\label{CondForTRP:3} %
$\Ld (m, \infty) = 0$ for all $m.$
\item\label{CondForTRP:4} %
$\Cs{2}{D}{\af}$ has a unique tracial state.
\item\label{CondForTRP:5} %
${\widehat{\af}}$ is trivial on $T (\Cs{2}{D}{\af}).$
\item\label{CondForTRP:6} %
If $\pi$ is the Gelfand-Naimark-Segal representation
associated with the unique tracial state $\ta$ on $D,$
then the automorphism ${\overline{\af}}$ of $\pi (D)''$
induced by $\af$ is outer.
\end{enumerate}
Moreover, if these conditions fail,
then $\Cs{2}{D}{\af}$ has exactly two extreme tracial states,
which are interchanged by ${\widehat{\af}}.$
\end{prp}

\begin{proof}
That (\ref{CondForTRP:1}) implies~(\ref{CondForTRP:2})
is Theorem~\ref{TRPDualToTAppRep}(\ref{TRPDualToTAppRep:1}).
To see that (\ref{CondForTRP:2}) implies~(\ref{CondForTRP:5}),
observe that ${\widehat{\af}}$ is tracially approximately
inner by Remark~\ref{TARImpTAI},
and apply Proposition~5.9 of~\cite{PhtRp1a}.
That (\ref{CondForTRP:4}) implies~(\ref{CondForTRP:5}) is obvious.
For the reverse,
it is evident from the description of the dual action on
the direct system of crossed products $C^* (\Zqt, D_n)$
in Lemma~\ref{Z2Struct} that $\Cs{2}{D}{\af}$ has only one
${\widehat{\af}}$-invariant tracial state.
The equivalence of~(\ref{CondForTRP:1}) and~(\ref{CondForTRP:6})
follows from Theorem~5.5 of~\cite{ELP}.

We now prove that~(\ref{CondForTRP:5}) implies~(\ref{CondForTRP:3}),
and also prove the last statement.
Assume~(\ref{CondForTRP:3}) fails.

Let ${\mathrm{tr}}_m$ denote the normalized trace on $M_m.$
Tracial states on $C^* (\Zqt, \, D, \, \gm )$ are in one to one
correspondence with sequences $(\ta_n)_{n \in \N}$ of tracial states
$\ta_n$ on $M_{t (n)} \oplus M_{t (n)}$ satisfying the
compatibility conditions $\ta_n \circ \ps_n = \ta_{n - 1}$
for all $n.$
The tracial state $\ta_n$ has the form
$\ta_n (a, b)
  = r_n {\mathrm{tr}}_{t (n)} (a) + s_n {\mathrm{tr}}_{t (n)} (b)$
for $r_n, s_n \in [0, 1]$ with $r_n + s_n = 1.$
To describe the compatibility condition,
for $\ld \in \R$ define the matrix
\[
T (\ld) = \frac{1}{2} \left( \begin{array}{cc}
           1 + \ld & 1 - \ld \\ 1 - \ld & 1 + \ld \end{array} \right)
        = \left( \begin{array}{cc}
               1  & 1  \\ 1  & -1 \end{array} \right)^{-1}
          \left( \begin{array}{cc}
               1  & 0  \\ 0  & \ld \end{array} \right)
          \left( \begin{array}{cc}
               1  & 1  \\ 1  & -1 \end{array} \right).
\]
{}From the second expression, we see that
$T (\ld \mu) = T (\ld) T (\mu)$ for $\ld, \, \mu \in \R,$
and that $T (1) = 1.$
{}From the first expression, we see that $k (n) T (\ld_n)$
is exactly the matrix $T_n$ of partial embedding multiplicities
of $\ps_n,$ as in Lemma~\ref{Z2Struct}(\ref{Z2Struct:2}).
Thus, the compatibility condition is exactly
$T (\ld_n) (r_n, s_n) = (r_{n - 1}, \, s_{n - 1})$
in $\R^2.$

For use below, we also observe that if $\ld \in [0, 1],$
if $(r_0, s_0), \, (r, s) \in \R^2$
satisfy $T (\ld) (r_0, s_0) = (r, s),$
and if $r_0, s_0 \in [0, 1]$ with $r_0 + s_0 = 1,$
then $r, s \in [0, 1]$ with $r + s = 1.$

We claim that a sequence $(r_n, s_n)_{n \in \N}$ corresponds to a
tracial state on the crossed product
$C^* (\Zqt, D, \gm )$ \ifo\  there is
$r \in [0, 1]$ such that for all $n \in \N$
we have $(r_n, s_n) = T (\Ld (n, \I)) (r, \, 1 - r).$
One direction is easy: given $r,$ the sequence $(r_n, s_n)_{n \in \N}$
defined by $T (\Ld (n, \I)) (r, \, 1 - r) = (r_n, s_n)$ clearly satisfies
$r_n, s_n \in [0, 1]$ and $r_n + s_n = 1,$
and the compatibility condition follows from the relation
$\Ld (n, \I) \ld_n = \Ld (n - 1, \, \I).$
For the other direction, choose $m \in \N$ such that
$\Ld (m, \infty) \neq 0.$
Define $r, s \in \R$ by
$(r, s) = T (\Ld (m, \infty))^{-1} (r_m, s_m).$
The compatibility condition guarantees that
$T (\Ld (n, \infty))^{-1} (r_n, s_n) = (r, s)$ for all $n \geq m.$
Now $\sum_{k = m + 1}^{\infty} \log (\ld_k)$ converges to
$\log (\Ld (m, \I)),$
so
\[
\log \left( \limi{n} \Ld (n, \infty) \right)
  = \limi{n} \sum_{l = n + 1}^{\infty} \log (\ld_l)
  = 0,
\]
whence $\limi{n} \Ld (n, \infty) = 1.$
Therefore $\limi{n} r_n = r$ and $\limi{n} s_n = s.$
It follows that $r, s \in [0, 1]$ and $r + s = 1.$
We already have $T (\Ld (n, \I)) (r, s) = (r_n, s_n)$
for $n \geq m,$
and the compatibility condition implies that this relation holds
for $n < m$ as well.
This completes the proof of the claim.

We now have an affine parametrization
of the tracial states on $\Cs{2}{D}{\gm}$ by $[0, 1],$
from which it is clear that
there are exactly two extreme tracial states.
That the dual action exchanges them
is clear from the identification of the dual action with the flip
in Lemma~\ref{Z2Struct}.
Thus, (\ref{CondForTRP:5})~fails.
We have also verified the last claim
whenever (\ref{CondForTRP:3})~fails.

We now prove that~(\ref{CondForTRP:3}) implies~(\ref{CondForTRP:1}).
Assuming~(\ref{CondForTRP:3}),
we verify the hypotheses of Theorem~\ref{ARPFromPosElts2}.
(Here, we use only a weak form of this theorem.)
It suffices to consider finite subsets
$F$ of the dense subalgebra
$\bigcup_{n = 1}^{\I} D_n \S D.$
Thus, let $F \S D_N,$ and let $\ep > 0.$
By hypothesis, we have $\Ld (N, \infty) = 0.$
So there exists $n$ such that
$\Ld (N, \, N + n) < \ep.$
Using Lemma~\ref{Condense},
write
\[
\bigotimes_{m = N + 1}^{N + n} \Ad ( p_m - q_m )
  = \Ad (p - q)
\]
for \pj s $p, q \in \bigotimes_{m = N + 1}^{N + n} M_{k (m)}$
such that $p + q = 1$
and $\rank (p) \geq \rank (q),$
and note that Lemma~\ref{Condense}(\ref{Condense:5}) implies
\[
\rank (p) - \rank (q) = \Ld (N, \, N + n) [\rank (p) + \rank (q)].
\]
Choose $v \in \bigotimes_{m = N + 1}^{N + n} M_{k (m)}$
such that $v v^* \leq p$ and $v^* v = q.$ 
Define $f_0, f_1 \in \bigotimes_{m = N + 1}^{N + n} M_{k (m)}$ by
\[
f_0 = \tfrac{1}{2} (v v^* + q + v + v^*)
\andeqn
f_1 = \tfrac{1}{2} (v v^* + q - v - v^*).
\]
Then one verifies that $f_0$ and $f_1$ are orthogonal \pj s
with sum $v v^* + q,$
and that moreover the unitary $w = p - q$ satisfies $w f_0 w^* = f_1$
and $w f_1 w^* = f_0.$
Now set
\[
e_0 = 1_{D_{N}} \otimes f_0,
\,\,  e_1 = 1_{D_{N}} \otimes f_1 \in D_{N + n} \S D.
\]
Then $e_0$ and $e_1$ commute exactly with every element of $D_N,$
and hence with every element of $F,$
and moreover $\af (e_0) = e_1$ and $\af (e_1) = e_0.$
Furthermore,
the unique tracial state $\ta$ satisfies
\[
\ta (1 - e_0 - e_1)
  = 1 - \frac{2 \cdot \rank (q)}{\rank (p) + \rank (q)}
  = \Ld (N, \, N + n)
  < \ep.
\]
So Theorem~\ref{ARPFromPosElts2} implies that $\af$ has the \tRp.
\end{proof}

\begin{prp}\label{CondForOut}
Let $\af \in \Aut (D)$ be a product type automorphism of order~$2$
as in Lemma~\ref{Z2Struct},
and adopt Notation~\ref{ProdTypeNtn}.
Then \tfae:
\begin{enumerate}
\item\label{CondForOut:1} %
The action of $\Zqt$ generated by $\af$ is outer, that is,
there is no unitary $u \in D$ such that $u^2 = 1$
and $u a u^* = \af (a)$ for all $a \in D.$
\item\label{CondForOut:2} %
The automorphism $\af$ is outer, that is,
there is no unitary $u \in D$ such that
$u a u^* = \af (a)$ for all $a \in D.$
\item\label{CondForOut:3} %
There are infinitely many $n \in \N$ such that $\ld_n < 1$
(that is, $q_n \neq 0$).
\item\label{CondForOut:4} %
$\Cs{2}{D}{\af}$ is simple.
\end{enumerate}
\end{prp}

\begin{proof}
We first prove that (\ref{CondForOut:1}) implies~(\ref{CondForOut:2}).
This uses only the fact that the center of $A$ is $\C \cdot 1.$
So suppose that $u \in D$ is a unitary such that
$u a u^* = \af (a)$ for all $a \in D.$
Putting $a = u,$ we get $\af (u) = u.$
Using $\af^2 = \id_A,$ we then find that $u^2$ is central.
Therefore $u^2 \in \C \cdot 1.$
So there is $\om \in \C$ such that $v = \om u$ satisfies $v^2 = 1,$
and clearly $v a v^* = \af (a)$ for all $a \in D.$

That (\ref{CondForOut:2}) implies~(\ref{CondForOut:4})
follows from Theorem~3.1 of~\cite{Ks1}.

If~(\ref{CondForOut:1}) fails,
it is well known that $\Cs{2}{D}{\af} \cong D \oplus D.$
Since $D \oplus D$ is not simple,
we have proved that (\ref{CondForOut:4}) implies~(\ref{CondForOut:1}).

It is easy to see that (\ref{CondForOut:1})
implies~(\ref{CondForOut:3}):
if (\ref{CondForOut:3})~fails,
then there is $n$ such that $q_m = 0$ for all $m > n,$
and
\[
u =
 (p_1 - q_1) \otimes (p_2 - q_2) \otimes \cdots
     \otimes (p_n - q_n) \otimes 1
\]
is a unitary such that $u^2 = 1$
and $u a u^* = \af (a)$ for all $a \in D.$

We prove that (\ref{CondForOut:3}) implies~(\ref{CondForOut:4}).
Using Lemma~\ref{Condense}(\ref{Condense:4}),
we may pass to a subsystem
and assume that $\rank (q_n) \neq 0$ for all $n.$
In this case, Lemma~\ref{Z2Struct}
shows that the partial embedding multiplicities
in the direct system of crossed products are all nonzero,
whence $\Cs{2}{D}{\af}$ is simple.
\end{proof}

We are now ready to give examples.
The first example is the most regular possible,
and is given to contrast with the remaining ones.

\begin{exa}\label{CAR1}
Let $\af$ be the automorphism of order~$2$ given by
\[
D = \bigotimes_{n = 1}^{\I} M_2 \andeqn
\af = \bigotimes_{n = 1}^{\I}
 \Ad \left( \begin{array}{cc} 1 & 0 \\ 0 &  -1 \end{array} \right).
\]
Then the action of $\Zqt$ generated by $\af$
is \sar\  and has the \sRp,
the crossed product is again the $2^{\infty}$~UHF algebra,
and the dual action is just another copy of the given action.
All this is easily proved using Proposition~\ref{CondForRP}
(and Lemma~\ref{Z2Struct} to precisely identify the
crossed product and dual action),
and is also a special case of Example~3.2 of~\cite{Iz}.
\end{exa}

\begin{exa}\label{CAR2}
Let $\af$ be the automorphism of order~$2$ given by
\[
D = \bigotimes_{n = 1}^{\I} M_{2^n} \andeqn
\af = \bigotimes_{n = 1}^{\I}
 \Ad \big( 1_{2^{n - 1} + 1} \oplus (- 1_{2^{n - 1} - 1}) \big).
\]
The automorphism in the $n$-th tensor factor is conjugation by a
diagonal unitary in which $2^{n - 1} + 1$ diagonal entries are equal
to $1$ and $2^{n - 1} - 1$ diagonal entries are equal to $- 1.$
Note that $D$ is the $2^{\infty}$~UHF algebra.

Following Notation~\ref{ProdTypeNtn},
we have $\ld_n \neq 0$ for all $n,$
but $\Ld (n, \infty) = 0$ for all $n.$
Applying Lemma~\ref{Z2Struct}
and Propositions~\ref{CondForRP} and~\ref{CondForTRP},
we obtain the following:
\begin{enumerate}
\item\label{CAR2:1} %
The action $\af$ has the \tRp\  but not the \sRp.
\item\label{CAR2:2} %
The action $\af$ is \sar.
\item\label{CAR2:3} %
The dual action ${\widehat{\af}}$ has the \sRp.
\item\label{CAR2:4} %
The dual action ${\widehat{\af}}$ is \tar\  but not \sar.
\item\label{CAR2:5} %
The dual action ${\widehat{\af}}$
is nontrivial on $K_0 (\Cs{2}{D}{\af}).$
\item\label{CAR2:6} %
$D$ is a UHF~algebra.
\item\label{CAR2:7} %
$\Cs{2}{D}{\af}$ is a simple AF~algebra
with unique tracial state but is not a UHF~algebra.
\end{enumerate}
\end{exa}

The following example was suggested by Izumi.
Some of the properties given here are folklore,
but we have been unable to find a reference for the proofs.
This example is a special case of an example used for other purposes
in Example~3.14 of~\cite{Iz0}.

\begin{exa}\label{CAR3}
Let $\af$ be the automorphism of order~$2$ given by
\[
D = \bigotimes_{n = 1}^{\I} M_{2^n} \andeqn
\af = \bigotimes_{n = 1}^{\I}
 \Ad ( 1_{2^{n} - 1} \oplus (- 1) ).
\]
The automorphism in the $n$-th tensor factor is conjugation by a
diagonal unitary in which $2^{n} - 1$ diagonal entries are equal
to $1$ and one diagonal entry is equal to $- 1.$
Note that $D$ is the $2^{\infty}$~UHF algebra.

Following Notation~\ref{ProdTypeNtn},
we have $\Ld (n, \infty) \neq 0$ for all $n.$
Applying Lemma~\ref{Z2Struct}
and Propositions~\ref{CondForRP}, \ref{CondForTRP}, and~\ref{CondForOut},
we obtain the following:
\begin{enumerate}
\item\label{CAR3:1} %
The action $\af$ is outer but does not have the \tRp.
\item\label{CAR3:2} %
The action $\af$ is \sar.
\item\label{CAR3:3} %
The dual action ${\widehat{\af}}$ has the \sRp.
\item\label{CAR3:4} %
The dual action ${\widehat{\af}}$ is not \tar.
\item\label{CAR3:5} %
The dual action ${\widehat{\af}}$
is nontrivial on $T (\Cs{2}{D}{\af}).$
\item\label{CAR3:6} %
$D$ is a UHF~algebra.
\item\label{CAR3:7} %
$\Cs{2}{D}{\af}$ is a simple AF~algebra with more than one tracial state.
\item\label{CAR3:8} %
In the factor representation of $D$ associated to the trace,
$\af$ becomes inner.
\end{enumerate}
\end{exa}

\begin{exa}\label{NotCAR}
Let $B$ be the $2^{\infty}$~UHF algebra,
and let $D$ be a UHF~algebra such that $B \otimes D \not\cong D.$
Then there exists an action of $\Zqt$ on $D$ with the \tRp,
but there is no action of $\Zqt$ on $D$ with the \sRp.

For the first statement,
write
\[
D = M_m \otimes \bigotimes_{n = 2}^{\I} M_{2 k (n) + 1}
\]
with $m$ a power of~$2$ and with $k (n) \geq 1$ for all $n.$
Then take
\[
\af = \id_{M_m} \otimes \bigotimes_{n = 1}^{\I}
 \Ad ( 1_{k (n) + 1} \oplus (- 1_{k (n)}) ).
\]
Following Notation~\ref{ProdTypeNtn},
one has $\ld_n \leq \frac{1}{3}$ for all $n \geq 2.$
Thus, Condition~(\ref{CondForTRP:3}) of Proposition~\ref{CondForTRP}
is satisfied.

Now suppose $\af$ is any action of $\Zqt$ on $D.$
Theorem~3.5 of~\cite{Iz2} implies, in particular, that if the
action had the \sRp,
then $D$ would be isomorphic to its tensor product with
the $2^{\infty}$~UHF algebra, contrary to hypothesis.
\end{exa}

\section{Blackadar's example has the tracial Rokhlin
  property}\label{Sec:Blackadar}

This section is devoted to the following example.

\begin{exa}\label{CAR4}
Let $A$ be the $2^{\infty}$ UHF algebra,
and let $\af$ be the automorphism
constructed in Section~5 of~\cite{Bl0}.
It follows from Corollary 5.3.2 of~\cite{Bl0}
and Takai duality~\cite{Tk} that $\Cs{2}{A}{\af}$ is not an AF~algebra.
We prove in Proposition~\ref{BlAutHasARP} below
that $\af$ generates an action of $\Zqt$ with the \tRp,
and in Proposition~\ref{CAR4_4} below
that this action does not have the \sRp,
and that the the dual action
is \tar\  %
but its generator induces a nontrivial automorphism of $K_1.$
By construction, the action of $\Zqt$ generated by $\af$
is \sar\  %
in the sense of Definition~\ref{AppRep}.
(The construction is recalled below.)
It follows from Theorem~\ref{RPDualToAppRep}(\ref{RPDualToAppRep:2})
that the dual action on $\Cs{2}{A}{\af}$ has
the \sRp.
\end{exa}

We remark that the methods used to prove the \tRp\  in this example
seem likely to be more
typical of proofs that actions on AH~algebras have the \tRp\  than
the methods used for Proposition~2.10 of~\cite{PhtRp2}.

We begin with a convenient description
of the construction in~\cite{Bl0}, following Section~5 there.
We make the convention that in any block matrix decomposition,
all blocks are to be the same size,
and we write $1_{n}$ for the identity of $M_n.$
We take the identification of $M_m \otimes M_n$ with $M_{m n}$
to send $a \otimes e_{j, k}$
to the $n \times n$ block matrix with $m \times m$ blocks,
of which the $(j, k)$ block is $a$ and the rest are zero.
Thus, $a \otimes 1 = \diag (a, a, \dots, a).$
Also, we identify the circle $S^1$ with $\R / \Z,$
and write elements of $C (S^1)$ as functions on $[0, 1]$
whose values at $0$ and $1$ are equal.

Following Definition 3.1.1 of~\cite{Bl0},
we choose a standard twice around embedding
$\ph^+ \colon C (S^1) \to C (S^1, M_2),$
given by choosing a \ct\  unitary path $c \in C ([0, 1], \, M_2)$
with
\[
c (0) = 1 \andeqn
c (1)
 = \left( \begin{array}{cc} 0 & 1 \\ 1 & 0 \end{array} \right),
\]
and then setting
\[
\ph^+ (f) (t) = c (t)
  \left( \begin{array}{cc} f \left( \frac{1}{2} t \right) & 0
    \\ 0 & f \left( \frac{1}{2} (t + 1) \right) \end{array} \right)
   c (t)^*
\]
for $f \in C (S^1).$
Further let $\ph^- \colon C (S^1) \to C (S^1, M_2)$
be the standard $- 2$ times around embedding
$\ph^- (f) (t) = \ph^+ (f) (1 - t).$
Extend everything, using the same notation, to embeddings
of $C (S^1, M_m)$ in $C (S^1, M_{2 m}),$
by using $\ph^+ \otimes \id_{M_m},$ etc.

Following Section~5 of~\cite{Bl0},
set $A_n = C (S^1, \, M_{4^n})$
and, remembering our convention on block sizes,
define a unitary in $M_{4^n} \S A_n$ by $u_n = \diag (1, \, -1).$
Further define $\ps_n \colon A_n \to A_{n + 1}$ by
\[
\ps_n \left( \begin{array}{cc} f_{1, 1} & f_{1, 2}
    \\ f_{2, 1} & f_{2, 2} \end{array} \right)
 = \left( \begin{array}{cccc}
\ph^+ (f_{1, 1}) &  0     & \ph^+ (f_{1, 2}) &  0      \\
  0     & \ph^- (f_{2, 2}) &  0     & \ph^- (f_{2, 1}) \\
\ph^+ (f_{2, 1}) &  0     & \ph^+ (f_{2, 2}) &  0      \\
  0     & \ph^- (f_{1, 2}) &  0     & \ph^- (f_{1, 1})
\end{array} \right).
\]
Theorem 4.1.1, Proposition 5.1.1, and Proposition 5.1.2 of~\cite{Bl0}
show that the direct limit of the $A_n$ using the maps $\ps_n$ is
the $2^{\infty}$ UHF algebra $A,$
and that the automorphisms $\af_n = \Ad (u_n)$ of $A_n$ define an
automorphism $\af$ of $A$ of order two.
It follows from Takai duality~\cite{Tk} and
Propositions 5.2.2 and 5.4.1 of~\cite{Bl0}
that $\Cs{2}{A}{\af}$ is isomorphic to the tensor product
of $A$ and the $2^{\infty}$ Bunce-Deddens algebra.

Further let $\io_n \colon M_{4^n} \to A_n$ be the embedding of
matrices as constant functions,
and define $\sm_n \colon M_{4^n} \to M_{4^{n + 1}}$ by
\[
\io_n \left( \begin{array}{cc} a_{1, 1} & a_{1, 2}
    \\ a_{2, 1} & a_{2, 2} \end{array} \right)
 = \left( \begin{array}{cccc}
a_{1, 1} \otimes 1_2 &  0     & a_{1, 2} \otimes 1_2 &  0      \\
  0     & a_{2, 2} \otimes 1_2 &  0     & a_{2, 1} \otimes 1_2 \\
a_{2, 1} \otimes 1_2 &  0     & a_{2, 2} \otimes 1_2 &  0      \\
  0     & a_{1, 2} \otimes 1_2 &  0     & a_{1, 1} \otimes 1_2
\end{array} \right).
\]
It is a consequence of the next lemma that
$\io_{n + 1} \circ \sm_n = \ps_n \circ \io_n.$
Moreover, we get an automorphism $\mu_n$ of $M_{4^n}$
by defining
$\mu_n = \Ad (u_n),$ and $\io_n \circ \mu_n = \af_n \circ \io_n.$

\begin{lem}\label{PresConst}
Let the notation be as above.
Let $\ep \geq 0,$ let $f \in A_n,$ and let $a \in M_{4^n}.$
Suppose $f (t) = a$ for all $t \in [\ep, \, 1 - \ep].$
Then $\ps_n (f)(t) = \sm_n (a)$
for all $t \in [2 \ep, \, 1 - 2 \ep].$
In particular, $\io_{n + 1} \circ \sm_n = \ps_n \circ \io_n.$
\end{lem}

\begin{proof}
Write
\[
f = \left( \begin{array}{cc} f_{1, 1} & f_{1, 2}
    \\ f_{2, 1} & f_{2, 2} \end{array} \right)
\andeqn
a = \left( \begin{array}{cc} a_{1, 1} & a_{1, 2}
    \\ a_{2, 1} & a_{2, 2} \end{array} \right).
\]
Then for each $j$ and $k,$
we have $f_{j, k} (t) = a_{j, k}$ for all $t \in [\ep, \, 1 - \ep].$
For $t \in [2 \ep, \, 1 - 2 \ep]$ we therefore get
\[
\left( \begin{array}{cc} f_{j, k} \left( \frac{1}{2} t \right) & 0
    \\ 0 & f_{j, k} \left( \frac{1}{2} (t + 1) \right) \end{array} \right)
= \left( \begin{array}{cc} a_{j, k} & 0
    \\ 0 & a_{j, k} \end{array} \right),
\]
which commutes with $c (t).$
\end{proof}

In the next lemma we show, roughly, that whenever an element
$f \in A_n = C (S^1, M_{4^n})$ is unitarily equivalent
in $C ([0, 1], \, M_{4^n}),$ via invariant unitaries,
to a function with small variation over intervals of length $\dt,$
then $\ps_n (f)$ is unitarily equivalent
in $C ([0, 1], \, M_{4^{n + 1}}),$ again via invariant unitaries,
to a function with small variation over intervals of length $2 \dt.$

\begin{lem}\label{AlmostConst}
Let the notation be as above.
Let $t \mapsto x (t)$ be a unitary element of $C ([0, 1], \, M_{4^n})$
such that $\mu_n (x (t)) = x (t)$ for every $t \in [0, 1].$
Then there exists
a unitary element $t \mapsto y (t)$ of $C ([0, 1], \, M_{4^{n + 1}})$
such that $\mu_{n + 1} (y (t)) = y (t)$ for every $t \in [0, 1],$
with the property that whenever
$\ep > 0,$ $\dt > 0,$ and $f \in A_n$ satisfy
\[
\| x (s) f (s) x (s)^* - x (t) f (t) x (t)^* \| < \ep
\]
for all $s, \, t \in [0, 1]$ such that $| s - t | < \dt,$
then
\[
\| y (s) \ps_n (f) (s) y (s)^* - y (t) \ps_n (f) (t) y (t)^* \| < \ep
\]
for all $s, \, t \in [0, 1]$ such that $| s - t | < 2 \dt.$
\end{lem}

\begin{proof}
The equation $\mu_n (x (t)) = x (t)$ implies that we can write
$x (t) = x_1 (t) \oplus x_2 (t)$ for unitaries
$x_1, \, x_2 \in C ([0, 1], \, M_{2^{2 n - 1}}).$
For $j = 1, 2$ define
\[
y_j (t) =
  \left( \begin{array}{cc} x_j \left( \frac{1}{2} t \right) & 0
    \\ 0 & x_j \left( \frac{1}{2} (t + 1) \right) \end{array} \right)
   c (t)^*.
\]
Then define
\[
y (t) = \diag ( y_1 (t), \, y_2 (1 - t), \, y_2 (t), \, y_1 (1 - t) ).
\]
Evidently $y$ is a unitary in $C ([0, 1], \, M_{4^{n + 1}})$
and $\mu_{n + 1} (y (t)) = y (t)$ for every $t \in [0, 1].$

To verify the conclusions of the lemma, it will simplify the
notation to conjugate everything by the permutation matrix
\[
w = \left( \begin{array}{cccc}
  1     &  0     &  0     &  0      \\
  0     &  0     &  1     &  0      \\
  0     &  0     &  0     &  1      \\
  0     &  1     &  0     &  0
\end{array} \right).
\]
(This conjugation is also used in Section~5 of~\cite{Bl0}.)
Thus, let
\[
{\widetilde{\ps}} (f) = w \ps_n (f) w^*
 = \left( \begin{array}{cccc}
\ph^+ (f_{1, 1}) & \ph^+ (f_{1, 2}) &  0     &  0      \\
\ph^+ (f_{2, 1}) & \ph^+ (f_{2, 2}) &  0     &  0      \\
  0     &  0     & \ph^- (f_{1, 1}) & \ph^- (f_{1, 2}) \\
  0     &  0     & \ph^- (f_{2, 1}) & \ph^- (f_{2, 2})
\end{array} \right),
\]
let
\[
{\widetilde{u}} = w u_{n + 1} w^* = \diag (1, \, -1, \, -1, \, 1),
\]
let
\[
{\widetilde{y}} (t) = w y (t) w^* 
 = \diag ( y_1 (t), \, y_2 (t), \, y_1 (1 - t), \, y_2 (1 - t) ),
\]
and similarly define ${\widetilde{\mu}},$ etc.
Note that ${\widetilde{\io}} = \io_{n + 1}.$

Let $\ep > 0,$ $\dt > 0,$ and $f \in A_n$ be as in the hypotheses.
Define
\[
g (t)
 = \left( \begin{array}{cc}
y_1 (t) \ph^+ (f_{1, 1}) (t) y_1 (t)^*
             & y_1 (t) \ph^+ (f_{1, 2}) (t) y_2 (t)^*  \\
y_2 (t) \ph^+ (f_{2, 1}) (t) y_1 (t)^*
             & y_2 (t) \ph^+ (f_{2, 2}) (t) y_2 (t)^*
\end{array} \right),
\]
and note that
\[
{\widetilde{y}} (t) {\widetilde{\ps}} (f) {\widetilde{y}} (t)^*
  = \diag (g (t), \, g (1 - t)).
\]
Accordingly, it suffices to prove that if
\[
\| x (s) f (s) x (s)^* - x (t) f (t) x (t)^* \| < \ep
\]
for all $s, \, t \in [0, 1]$ such that $| s - t | < \dt,$
then
\[
\| g (s) - g (t) \| < \ep
\]
for all $s, \, t \in [0, 1]$ such that $| s - t | < 2 \dt.$

Let $v$ be the permutation matrix
\[
v = \left( \begin{array}{cccc}
  1     &  0     &  0     &  0      \\
  0     &  0     &  1     &  0      \\
  0     &  1     &  0     &  0      \\
  0     &  0     &  0     &  1
\end{array} \right).
\]
When one calculates $v g (t) v^*$ by
substituting the formulas for $y_j (t)$ and $\ph^+$ in the expression
for $g (t),$
the factors $c (t)$ and $c (t)^*$ all cancel out,
and the final answer is
\[
v g (t) v^* = \ts{
  \diag \left(
x \left( \frac{1}{2} t \right) f \left( \frac{1}{2} t \right)
             x \left( \frac{1}{2} t \right)^*, \,\,
x \left( \frac{1}{2} (t + 1) \right)
             f \left( \frac{1}{2} (t + 1) \right)
                       x \left( \frac{1}{2} (t + 1) \right)^* \right). }
\]
Since we are assuming
\[
\| x (s) f (s) x (s)^* - x (t) f (t) x (t)^* \| < \ep
\]
for all $s, \, t \in [0, 1]$ such that $| s - t | < \dt,$
it is immediate that
$| s - t | < 2 \dt$ implies
\[
\| v g (s) v^* - v g (t) v^* \| < \ep,
\]
whence also
\[
\| g (s) - g (t) \| < \ep,
\]
as desired.
\end{proof}

\begin{prp}\label{BlAutHasARP}
The automorphism $\af$ of Example~\ref{CAR4}
generates an action of $\Zqt$ with the \tRp.
\end{prp}

\begin{proof}
Let the notation be as before Lemma~\ref{PresConst}.
Let $\ta$ be the unique tracial state on $A = \dirlim A_n.$
Define a tracial state $\ta_n$ on $A_n$ by
\[
\ta_n (f) = \int_0^1 {\mathrm{tr}}_{4^n} (f (t ) ) \, d t,
\]
where ${\mathrm{tr}}_{m}$ is the normalized trace on $M_m.$
Then one checks that $\ta_{n + 1} \circ \ps_n = \ta_n$ for all $n.$
It follows from the uniqueness of $\ta$ that $\ta |_{A_n} = \ta_n$
for all $n.$

We use Theorem~\ref{ARPFromPosElts2} to verify the \tRp.
So let $F \S A$ be finite and let $\ep > 0.$
Choose $m$ and a finite set $S_0 \S A_m$ such that
every element of $F$ is within $\frac{1}{8} \ep$ of an element of $S_0.$

The set $S_0$ is a uniformly equicontinuous set of functions from
$[0, 1]$ to $M_{4^m},$ so there is $\dt > 0$ such that
whenever $s, \, t \in [0, 1]$ satisfy $| s - t | < \dt,$ then
\[
\| f (s) - f (t) \| < \ts{ \frac{1}{8} } \ep
\]
for all $t \in [0, 1]$ and all $f \in S_0.$
Choose $n \in \N$ with $n \geq m$ and so large that $2^{n - m} \dt > 1.$
Apply Lemma~\ref{AlmostConst} a total of $n - m$ times,
the first time with $x (t) = 1$ for all $t,$ 
obtaining after the last application
a \ct\  unitary path $t \mapsto z (t)$ in $C ([0, 1], \, M_{4^n})$
such that $\mu_{n} (z (t)) = z (t)$ for every $t \in [0, 1].$
Replacing $z (t)$ by $z (0)^* z (t),$
we may clearly assume that $z (0) = 1.$
Then, in particular,
\[
\| z (t)^* f (0) z (t) - f (t) \| < \ts{ \frac{1}{8} } \ep
\]
for all $t \in [0, 1]$ and all $f \in S_0.$
Recall that we identify $C (S^1, B)$ with the set of
functions $f \in C ([0, 1], \, B)$ such that $f (0) = f (1).$
Since the fixed point algebra $A_n^{\af_n} = C (S^1, M_{4^n})^{\af_n}$
is just $C (S^1, \, M_{4^n}^{\mu_n}),$
and since $M_{4^n}^{\mu_n}$ is \fd,
there is an $\af_n$-invariant unitary $y \in A_n$ such that
$y (t) = z (t)$ for $t \in \left[0, \, 1 - \frac{1}{8} \ep \right].$
Then for each $f \in S_0,$ regarded as a subset of $A_n,$ there exists
$g \in A_n$ such that $\| y^* g y - f \| < \ts{ \frac{1}{4} } \ep$
and $g (t) = f (0)$
for $t \in \left[0, \, 1 - \frac{1}{8} \ep \right].$
Let $S$ be the set of all elements $g$ obtained in this way from
elements of $S_0.$
In particular, for every $a \in F$ there is $g \in S$ such that
$\| a - y^* g y \| < \frac{1}{2} \ep.$

We claim that there are orthogonal positive elements
$b_0, \, b_1 \in A_{n + 1} \S A$ such that
$b_j g = g b_j$ for all $g \in S,$
and such that
\[
\af_{n + 1} (b_0) = b_1, \,\,\,\,\,\,
\af_{n + 1} (b_1) = b_0, \,\,\,\,\,\,
0 \leq b_0, \, b_1 \leq 1, \andeqn
0 \leq \ta (1 - b_0 - b_1) < \ep.
\]
For this purpose, it suffices to use in place of $\ps_n$ the
unitarily equivalent \hm\  %
\[
{\widetilde{\ps}} = w \ps_n (-) w^*
  \colon C (S^1, M_{4^n}) \to C (S^1, M_{4^{n + 1}})
\]
in the proof of Lemma~\ref{AlmostConst}
(called $\om_n$ in the proof of Proposition~5.1.1 of~\cite{Bl0}),
and to use in place of $\af_{n + 1}$ the automorphism
\[
{\widetilde{\af}} = \Ad (w) \circ \af_{n + 1} = \Ad ({\widetilde{u}})
 = \Ad (\diag (1, \, -1, \, -1, \, 1) ).
\]
Note that this change does not require any change in the formula
for the trace $\ta_{n + 1},$
and also does not affect the first part of
the conclusion of Lemma~\ref{PresConst}.
Accordingly, if
\[
g = \left( \begin{array}{cc} g_{1, 1} & g_{1, 2}
    \\ g_{2, 1} & g_{2, 2} \end{array} \right)
  \in S
\]
then ${\widetilde{\ps}} (g) \in  C (S^1, M_{4^{n + 1}})$ satisfies
\[
{\widetilde{\ps}} (g) (t) = \left( \begin{array}{cccc}
g_{1, 1} (0) & g_{1, 2} (0) &  0     &  0      \\
g_{2, 1} (0) & g_{2, 2} (0) &  0     &  0      \\
  0     &  0     & g_{1, 1} (0) & g_{1, 2} (0) \\
  0     &  0     & g_{2, 1} (0) & g_{2, 2} (0)
\end{array} \right)
 = \left( \begin{array}{cc} g (0) & 0 \\ 0 & g (0) \end{array} \right)
\]
for $t \in \left[\frac{1}{4} \ep, \, 1 - \frac{1}{4} \ep \right].$

Now set
\[
p_0 = \frac{1}{2} \left( \begin{array}{cccc}
  1     &  0     &  1     &  0      \\
  0     &  1     &  0     &  1      \\
  1     &  0     &  1     &  0      \\
  0     &  1     &  0     &  1
\end{array} \right)
\andeqn
p_1 = \frac{1}{2} \left( \begin{array}{cccc}
  1     &  0     & - 1    &  0      \\
  0     &  1     &  0     & - 1     \\
 - 1    &  0     &  1     &  0      \\
  0     & - 1    &  0     &  1
\end{array} \right),
\]
both in $M_4 (M_{4^n}).$
In $2 \times 2$ block form, we can write
\[
{\widetilde{u}}
 = \left( \begin{array}{cc} s & 0 \\ 0 & - s \end{array} \right)
\,\,\,\,\,\, {\mbox{with}} \,\,\,\,\,\,
s = \left( \begin{array}{cc} 1 & 0 \\ 0 & - 1 \end{array} \right),
\]
\[
p_0 = \frac{1}{2}
  \left( \begin{array}{cc} 1 & 1 \\ 1 & 1 \end{array} \right),
\andeqn
p_1 = \frac{1}{2}
  \left( \begin{array}{cc} 1 & - 1 \\ - 1 & 1 \end{array} \right).
\]
With these formulas, it is easy to check that
$p_0$ and $p_1$ are \pj s with $p_0 + p_1 = 1,$
that $p_0$ and $p_1$ commute with ${\widetilde{\ps}} (g) (t)$
for every $g \in S$
and $t \in \left[ \frac{1}{4} \ep, \, 1 - \frac{1}{4} \ep \right],$
and that $\Ad ({\widetilde{u}})$ exchanges $p_0$ and $p_1.$

Now choose and fix a \cfn\  $h \colon [0, 1] \to [0, 1]$
such that $h (t) = 0$ for
$t \not\in \left[ \frac{1}{4} \ep, \, 1 - \frac{1}{4} \ep \right]$
and $h (t) = 1$ for
$t \in \left[ \frac{1}{2} \ep, \, 1 - \frac{1}{2} \ep \right],$
and define $b_j (t) = h (t) p_j$ for $j = 0, \, 1.$
Then $b_0$ and $b_1$ are positive elements with $b_0, \, b_1 \leq 1,$
which commute with ${\widetilde{\ps}} (g)$ for every $g \in S,$
which satisfy $b_0 b_1 = 0,$
such that $\Ad ({\widetilde{u}})$ exchanges $b_0$ and $b_1,$
and such that $0 \leq \ta_{n + 1} (1 - b_0 - b_1) < \ep.$
This proves the claim above.

We return to the use of $\ps_{n + 1},$
and we let $b_0, \, b_1 \in A_{n + 1} \S A$ be as in the claim
(rather than its proof).
In $A_{n + 1} \S A,$ define $a_0 = y^* b_0 y$ and $a_1 = y^* b_1 y.$
Since $\af (y) = y,$ it follows that
$a_0, \, a_1 \in A_{n + 1} \S A$ satisfy
$a_j y^* g y = y^* g y a_j$ for all $g \in S,$
and
\[
a_0 a_1 = 0, \,\,\,\,\,\,
\af_{n + 1} (a_0) = a_1, \,\,\,\,\,\,
0 \leq a_0, \, a_1 \leq 1, \andeqn
0 \leq \ta (1 - a_0 - a_1) < \ep.
\]
For $a \in F$ choose $g \in S$ such that
$\| a - y^* g y \| < \frac{1}{4} \ep.$
Then
\[
\| [a_j, a ] \| \leq 2 \| a - y^* g y \| + \| [a_j, \, y^* g y] \|
  < 2 \left( \ts{\frac{1}{2} \ep} \right) + 0 = \ep.
\]
This completes the verification of the hypotheses of
Theorem~\ref{ARPFromPosElts2}, so it follows that $\af$ has the \tRp.
\end{proof}

\begin{prp}\label{CAR4_4}
Let $\af \in \Aut (A)$ be as in Example~\ref{CAR4}.
Then:
\begin{enumerate}
\item\label{CAR4_4:1} %
The action of $\Zqt$ generated by $\af$ does not have the \sRp.
\item\label{CAR4_4:2} %
The dual action on $\Cs{2}{A}{\af}$ has the \sRp.
\item\label{CAR4_4:3} %
The dual action is \tar.
\item\label{CAR4_4:4} %
The generating automorphism ${\widehat{\af}}$ of the dual action
acts nontrivially on $K_1 (\Cs{2}{A}{\af}).$
\end{enumerate}
\end{prp}

\begin{proof}
We have already observed in
Example~\ref{CAR4} that $\Cs{2}{A}{\af}$ is not AF.
Therefore~(\ref{CAR4_4:1}) follows from Theorem~2.2 of~\cite{PhtRp1a}.

It is immediate from the discussion following
Example~\ref{CAR4} that
the action of $\Zqt$ generated by $\af$
is \sar\  %
in the sense of Definition~\ref{AppRep}.
Part~(\ref{CAR4_4:2}) therefore follows from
Theorem~\ref{RPDualToAppRep}(\ref{RPDualToAppRep:2}).

Part~(\ref{CAR4_4:3}) follows from
Proposition~\ref{BlAutHasARP} and
Theorem~\ref{TRPDualToTAppRep}(\ref{RPDualToAppRep:1}).

It remains to prove~(\ref{CAR4_4:4}).
We continue to follow the notation introduced after
Example~\ref{CAR4}.
Let $B_n$ be the fixed point algebra
\[
A_n^{\af_n}
 = C (S^1, \, M_{2^{2 n - 1}}) \oplus C (S^1, \, M_{2^{2 n - 1}})
 \S C (S^1, \, M_{4^{n}}),
\]
with the embedding being as $2 \times 2$ block diagonal matrices.
Let $B = \dirlim B_n,$ which is also equal to $A^{\af}.$
Let $\bt_n \in \Aut (B_n)$ be $\bt_n (f, g) = (g, f).$
By Proposition~5.2.2 of~\cite{Bl0} and the preceding discussion, 
there is a corresponding automorphism $\bt$ of the direct limit,
$A \cong \Cs{2}{B}{\bt},$
and the isomorphism can be chosen so that $\af$ generates the
dual action.
By Takai duality~\cite{Tk}, it therefore suffices to show that
$\bt$ is nontrivial on $K_1 (B).$

Following the discussion after Corollary~5.3.2 of~\cite{Bl0},
let $v \in B_1 \S A_1 = C (S^1, M_4)$ be the unitary
\[
v (t) = \diag \left( e^{2 \pi i t}, \, e^{2 \pi i t},
   \, e^{- 2 \pi i t}, \, e^{- 2 \pi i t} \right).
\]
As there, the image of $[v]$ in $K_1 (B)$ is nonzero.
It follows from Proposition~5.3.1 of~\cite{Bl0}
that $K_1 (B)$ is torsion free,
and one checks that $[\bt_1 (v)] = - [v],$
so $\bt_* ([v]) = - [v] \neq [v].$
\end{proof}

\section{Torsion in $K_0$ of the crossed product}\label{Sec:Tor}

\indent
In this section, we give two examples related to 10.11.3
of~\cite{Bl}.
Let $A$ be the simple unital AF~algebra with
$K_0 (A) \cong \Z \big[ \frac{1}{2} \big] \oplus \Z,$ with
$\Z \big[ \frac{1}{2} \big]$ given the order from $\R,$
with the strict
order from the first coordinate on the direct sum, and with
$[ 1_A ] = (1, 0).$
A proof is outlined in 10.11.3 of~\cite{Bl} that if there is an
automorphism $\af$ of $A$ of order~$2$ such that
$\af_* ( \et, k) = ( \et, \, -k)$ on $K_0, (A),$ then
$C^* ( \Z_2, A, \af)$ is not AF because
$K_1 ( C^* ( \Z_2, A, \af )) \neq 0$ or $K_0 ( C^* ( \Z_2, A, \af ))$
has torsion.
Using the K-theory of any odd UHF~algebra in place of
$\Z \big[ \frac{1}{2} \big],$
we give examples of actions of $\Z_2$
showing that both possibilities can occur,
and which moreover have the \tRp.
To our knowledge, no previous example was known of an action of a
finite group on an AF~algebra such that the K-theory of the crossed
product has torsion.
The fact that the actions have the \tRp\   allows
us to show that the crossed products are classifiable.

\begin{exa}\label{E:Torsion}
Let $m \in \N.$
Define $h \colon S^{2 m} \to S^{2 m}$ by $h (x) = -x,$ and let
$\bt \in \Aut (C (S^{2 m} ))$ be the corresponding automorphism of
order~$2.$
For $r \in \N$ and $b \in S^{2 m},$ define
$\ps_{r, b} \colon C (S^{2 m}) \to M_{2 r + 1} \otimes C (S^{2 m})$
by
\[
\ps_{r, b} (f) (x)
 = \diag \big( f (x), \, f (b), \, f (h (b)), \, f (b), \, f (h (b)),
            \, \ldots, \, f (b), \, f (h (b)) \big)
\]
for $x \in S^{2 m},$ where $f (b)$ and $f (h (b))$ each occur $r$ times.
Choose a dense sequence $(x (n))_{n \in \N}$ in $S^{2 m},$ and choose a
sequence $(r (n))_{n \in \N}$ of positive integers such that
$\lim_{n \to \I} r (n) = \I.$
Set
\[
s (n) = [2 r (1) + 1] [2 r (2) + 1] \cdots [2 r (n) + 1],
\]
and set
$A_n = M_{s (n)} \otimes C (S^{2 m}),$
which, when appropriate, we think of as
\[
M_{2 r (1) + 1} \otimes M_{2 r (2) + 1}
        \otimes \cdots \otimes M_{2 r (n) + 1}
        \otimes C (S^{2 m}).
\]
Define $\ph_n \colon A_{n - 1} \to A_n$ by
$\ph_n = \id_{M_{s (n - 1)}} \otimes \ps_{x (n), \, r (n)}.$
Then set $A = \dirlim A_n.$

For $r \in \N$ define a unitary $w_r \in M_{2 r + 1}$ by
\[
w_r = \diag \left( 1, \,\,
 \left( \begin{array}{cc} 0 & 1 \\ 1 & 0 \end{array} \right), \,\,
 \left( \begin{array}{cc} 0 & 1 \\ 1 & 0 \end{array} \right), \,\,
 \ldots, \,\,
 \left( \begin{array}{cc} 0 & 1 \\ 1 & 0 \end{array} \right)
 \right).
\]
Then define an automorphism $\af_n \in \Aut (A_n)$ of order~$2$ by
\[
\af_n
 = \Ad (w_{r (1)} \otimes w_{r (2)} \otimes \cdots \otimes w_{r (n)} )
\otimes \bt.
\]
One checks that
$\ph_n \circ \af_{n - 1} = \af_n \circ \ph_n,$
so that the automorphisms $\af_n$
define an automorphism $\af \in \Aut (A)$ of order~$2.$
\end{exa}

\begin{prp}\label{P:Torsion}
Let $A$ and $\af \in \Aut (A)$ be as in Example~\ref{E:Torsion}.
Then:
\begin{enumerate}
\item\label{P:Torsion:1} %
$A$ is a simple unital AF~algebra with a unique tracial state.
\item\label{P:Torsion:2} %
$\Cs{2}{A}{\af}$ is a simple unital AH~algebra with no dimension
growth, tracial rank zero, and a unique tracial state.
\item\label{P:Torsion:3} %
The action of $\Zqt$ generated by $\af$
has the \tRp\  but not the \sRp.
\item\label{P:Torsion:4} %
The action of $\Zqt$ generated by $\af$
is \tar\  but not \sar.
\item\label{P:Torsion:5} %
There is a dense subgroup $G$ of $\R,$ contained in $\Q,$
such that $K_0 (A) \cong G \oplus \Z,$
with the strict order from the first coordinate,
and $\af_* (\et, k) = (\et, \, - k)$ for $\et \in G$ and $k \in \Z.$
\item\label{P:Torsion:6} %
The torsion subgroup of $K_0 (\Cs{2}{A}{\af})$
is isomorphic to $\Z_{2^m},$
and we have $K_1 (\Cs{2}{A}{\af}) = 0.$
\end{enumerate}
\end{prp}

\begin{proof}
We continue using the notation of Example~\ref{E:Torsion}.
The C*-algebra $A$ is almost of the form considered in~\cite{Gd},
the only difference being that each map in the direct system uses point
evaluations at two different points rather than just one.
However, most of the proofs are still valid.
In particular, $A$ is simple (as in Lemma~1 of~\cite{Gd}), and $A$ has
real rank zero (by the same proof as for Theorem~9 of~\cite{Gd}).
In addition,
the K-theory computation of Theorem~13 of~\cite{Gd} is still valid.
We conclude that $K_1 (A) = 0.$
Moreover, if we let $D$ be the UHF~algebra
$D = \bigotimes_{n = 1}^{\I} M_{2 r (n) + 1},$
and use the fact that the
reduced $K_0$-group ${\widetilde{K}}^0 (S^{2 m})$ is $\Z,$ we find that
$K_0 (A) \cong K_0 (D) \oplus \Z.$
The order unit is $( [1], 0),$ and the order is the strict order from
the first coordinate.
That is,
\[
K_0 (A)_{+}
  = \{ (\et, k ) \in K_0 (D) \oplus \Z \colon \et > 0 \}
          \cup \{ (0, 0) \}.
\]
Lin's classification theorem (Theorem~5.2 of~\cite{Ln15}),
in the form given in Proposition~3.7 of~\cite{PhtRp2},
therefore applies to this direct limit,
and we conclude that $A$ is an AF~algebra.
Note that it has a unique tracial state $\ta.$
Since $\bt$ is multiplication by $-1$ on
${\widetilde{K}}^0 (S^{2 m})$
(it reverses the sign of the Bott element),
it follows that, with $K_0 (A)$ as described above,
$\af_* (\et, k) = (\et, \, - k)$ for $\et \in K_0 (D)$ and $k \in \Z.$
We have proved~(\ref{P:Torsion:1}) and~(\ref{P:Torsion:5}).

We show that the action of $\Z_2$ generated by $\af$ has the \tRp,
which is the first part of~(\ref{P:Torsion:3}).
By Lemma~\ref{ARPFromPosElts},
it suffices to show that for every $n$ and every $\ep > 0,$ there is a
projection $e \in A$ such that $e \af (e) = 0,$ $e$ and $\af (e)$
commute with every element of $A_n,$ and $\ta (1 - e - \af (e)) < \ep.$
Choose $k \geq n + 1$ such that $1 / [2 r (k) + 1] < \ep,$ take
\[
e_0 = \diag \left( 0, \, \,
\left( \begin{array}{cc} 1 & 0 \\ 0 & 0 \end{array} \right), \, \,
\left( \begin{array}{cc} 1 & 0 \\ 0 & 0 \end{array} \right), \, \,
\ldots, \,\,
\left( \begin{array}{cc} 1 & 0 \\ 0 & 0 \end{array} \right)
 \right) \in M_{2 r (k) + 1},
\]
and take
\[
e = 1_{M_{s (k - 1)}} \otimes e_0 \otimes 1_{C (S^{2 m})} \in A_k.
\]

We now compute the crossed product.
Since the action of $\Z_2$ on $S^{2 m}$ is free, and
$S^{2 m} / \Z_2$ is the real projective space $\R P^{2 m},$
there is a strong Morita equivalence
between the fixed point algebra $C( S^{2 m})^{\Z_2} = C (\R P^{2 m})$
and the crossed product $B = C^* ( \Z_2, S^{2 m}, h).$
(For example, see Situation~2 in~\cite{Rf}.
Our case is exactly the specific example described there.)
So $C^* (\Z_2, A, \af)$
is an AH~algebra with no dimension growth,
which is the first part of~(\ref{P:Torsion:2}).
Tracial rank zero follows from Theorem~2.6 of~\cite{PhtRp1a},
and simplicity follows from Corollary~1.6 of~\cite{PhtRp1a}.
Because $A$ has a unique tracial state,
Proposition~5.7 of~\cite{ELP} implies that
$C^* (\Z_2, A, \af)$ has a unique tracial state.

We have $K_* (B) \cong K^* (\R P^{2 m}).$
Thus (see Proposition 2.7.7 of~\cite{At}),
$K_1 (B) = 0$ and $K_0 (B) \cong \Z \oplus \Z_{2^m}.$
Moreover, for any unital
homomorphism $\ep \colon B \to M_2,$ the isomorphism identifies
$\ep_* \colon \Z \oplus \Z_{2^m} \to \Z$ with the projection to the
first factor.

Fix $n,$ and set
$y = w_{r (1)} \otimes w_{r (2)} \otimes \cdots \otimes w_{r (n - 1)}.$
Then $\af_{n - 1} = \Ad (y) \otimes \bt$ and
$\af_n = \Ad (y) \otimes \Ad (w_r) \otimes \bt.$
The actions of $\Z_2$ these automorphisms generate are exterior
equivalent to those generated by $\id_{M_{s (n - 1)}} \otimes \bt$
and $\id_{M_{s (n - 1)}} \otimes \Ad (w_r) \otimes \bt,$
in a manner which respects $\ph_n.$
Since crossed products by exterior equivalent actions are naturally
isomorphic (see, for example, the proof of Theorem 2.8.3(5)
of~\cite{Ph}), we may identify the crossed product map
${\overline{\ph}}_n \colon B_{n - 1} \to B_n$ with the tensor product
$\id_{M_{s (n - 1)}} \otimes {\overline{\ps}}_{r (n), \, x (n)}$ of
$\id_{M_{s (n - 1)}}$ and the map
\[
{\overline{\ps}}_{r (n), \, x (n)}
 \colon B \to
  C^* ( \Z_2, \,\, M_{2 r (n) + 1} \otimes C( S^{2 m}),
                            \,\, \Ad (w_r) \otimes \bt)
\]
on crossed products induced by $\ps_{r (n), \, x (n)}.$
The codomain of ${\overline{\ps}}_{r (n), \, x (n)}$
is $M_{2 r (n) + 1} \otimes B$
by exterior equivalence.
Set
\[
w = \left( \begin{array}{cc} 0 & 1 \\ 1 & 0 \end{array} \right).
\]
Let $\ep \colon C ( S^{2 m} ) \to M_2 ( \C )$ be
$\ep (f) = \diag ( f (x (n)), \, f (h (x (n))),$ let
$\io \colon M_2 \to M_2 \otimes C ( S^{2 m})$ be
$a \mapsto a \otimes 1,$ let $\Z_2$ act on $M_2$ via $\Ad (w),$ and
let ${\overline{\ep}}$ and ${\overline{\io}}$ be the induced maps on
crossed products.
Using naturality,
we can identify ${\overline{\ps}}_{r (n), \, x (n)}$ with
\[
a \mapsto
 \diag (a, \, {\overline{\io}} \circ {\overline{\ep}} (a), \,
      {\overline{\io}} \circ {\overline{\ep}} (a), \, \ldots, \,
      {\overline{\io}} \circ {\overline{\ep}} (a) ),
\]
with ${\overline{\io}} \circ {\overline{\ep}} (a)$ occurring $r (n)$
times.
We have $C^* ( \Z_2, \, M_2, \, \Ad (w)) \cong M_2 \oplus M_2,$
and the map ${\overline{\ep}}$ is unital.
Identify
$K_0 ( M_2 \oplus M_2 )$ with $\Z \oplus \Z$ in the obvious way.
Identify $K_0 (B)$ with $\Z \oplus \Z_{2^m}$ as above,
and recall our analysis of the K-theory of unital \hm s $B \to M_2.$
Then we have ${\overline{\ep}}_* (k, l) = (k, k)$
for $k \in \Z$ and $l \in \Z_{2^m}.$
Putting the summands together, we find that
$\big( {\overline{\ps}}_{r (n), \, x (n)} \big)_* (k, l)
  = \big( ( 2 r (n) + 1) k, \, l \big).$
It follows that, with $D$ being the same UHF~algebra as above, we have
\[
K_0 (C^* (\Z_2, A, \af)) \cong K_0 (D) \oplus \Z_{2^m}
\andeqn
K_1 ( C^* (\Z_2, A, \af )) = 0.
\]
We have proved~(\ref{P:Torsion:6}).
Since $C^* (\Z_2, A, \af )$ is now clearly not~AF,
an application of Theorem~2.2 of~\cite{PhtRp1a}
show that the action of $\Z_2$ generated by $\af$
does not have the \sRp,
which completes the proof of~(\ref{P:Torsion:3}).

It remains to prove~(\ref{P:Torsion:4}).
We claim that the map $K_0 (A) \to K_0 (C^* (\Z_2, A, \af))$
is not injective.
It is clear that all elements $(\et, 0) \in K_0 (D) \oplus \Z_{2^m},$
with $\et > 0,$ are positive in $K_0 (C^* (\Z_2, A, \af)).$
The infinitesimals of $K_0 (C^* (\Z_2, A, \af))$ are therefore
contained in $0 \oplus \Z_{2^m}.$
Since infinitesimals map to infinitesimals,
and the summand $0 \oplus \Z \S K_0 (D) \oplus \Z \cong K_0 (A)$
consists of infinitesimals,
the claim follows.
Since $A$ is the fixed point algebra of $C^* (\Z_2, A, \af)$
under the dual action ${\widehat{\af}},$
Theorem~3.13 of~\cite{Iz}
implies that ${\widehat{\af}}$ does not have the \sRp.
So, by Theorem~\ref{RPDualToAppRep}(\ref{RPDualToAppRep:2}),
the action generated by $\af$ is not \sar.

On the other hand,
$\af$ is strongly tracially approximately inner
by Theorem~6.6 of~\cite{PhtRp1a},
and we know $\af$ generates an action with the \tRp,
so Theorem~4.6 of~\cite{PhtRp1a} implies that this action is \tar.
\end{proof}

Except for the fact that we have the K-theory of an odd UHF~algebra in
place of $\Z \big[ \frac{1}{2} \big],$
we have constructed an automorphism
of order~$2$ of the type suggested in 10.11.3 of~\cite{Bl}.

\begin{rmk}\label{R:Torsion_Order}
We continue with the notation of Example~\ref{E:Torsion}.
Although we do not carry out the details, it is not hard to show that
the order on $K_0 (C^* (\Z_2, A, \af)) \cong K_0 (D) \oplus \Z_{2^m}$
is the strict order from the first coordinate.
We saw in the proof of Proposition~\ref{P:Torsion}
that Lin's classification theorem
(Theorem~5.2 of~\cite{Ln15}) applies.
(We could also use the classification theorem of~\cite{EGL},
without appealing to the \tRp.)
It follows, for example, that $C^* ( \Z_2, A, \af)$ is
stably isomorphic to
the direct limit of a system of the form
\[
C (\R P^{2 m}) \longrightarrow M_{2 r (1) + 1} \otimes C ( \R P^{2 m})
 \longrightarrow M_{(2 r (1) + 1)(2 r (2) + 1)} \otimes C ( \R P^{2 m})
 \longrightarrow \cdots,
\]
using one copy of the identity map and $2 r (n)$ point evaluations at
the $n$-th stage.
\end{rmk}

\begin{rmk}\label{R:Torsion_Dual}
It follows from Theorems \ref{RPDualToAppRep} and~\ref{TRPDualToTAppRep}
that the dual of the action in Example~\ref{E:Torsion} also
has the \tRp\  but not the \sRp,
and is \tar\  but not \sar.
\end{rmk}

\begin{exa}\label{E:NoTor}
This example is the same as Example~\ref{E:Torsion},
except that we use a different homeomorphism $h,$
namely
$h (x_0, x_1, \ldots, x_{2 m}) = (-x_0, x_1, \ldots, x_{2 m}).$
We choose our dense
sequence to be disjoint from the fixed points of $h$;
this is possible because the complement of the fixed point set is a
dense open set.
We define the maps $\ps_{r, b},$ $\ph_{n, \bt},$ and $\af_n,$ and the
UHF~algebra $D,$ by the same formulas as in Example~\ref{E:Torsion}.
Let $A$ be the direct limit, as before,
and let $\af$ be the automorphism of order~$2$
obtained as the direct limit of the automorphisms~$\af_n.$
\end{exa}

\begin{prp}\label{P:NoTor}
Let $A$ and $\af \in \Aut (A)$ be as in Example~\ref{E:NoTor}.
Then:
\begin{enumerate}
\item\label{P:NoTor:1} %
$A$ is a simple unital AF~algebra with a unique tracial state.
\item\label{P:NoTor:2} %
$\Cs{2}{A}{\af}$ is a simple unital AH~algebra with no dimension
growth, tracial rank zero, and a unique tracial state.
\item\label{P:NoTor:3} %
The action of $\Zqt$ generated by $\af$
has the \tRp\  but not the \sRp.
\item\label{P:NoTor:4} %
The action of $\Zqt$ generated by $\af$
is \tar\  but not \sar.
\item\label{P:NoTor:5} %
There is a dense subgroup $G$ of $\R,$ contained in $\Q,$
such that $K_0 (A) \cong G \oplus \Z,$
with the strict order from the first coordinate,
and $\af_* (\et, k) = (\et, \, - k)$ for $\et \in G$ and $k \in \Z.$
\item\label{P:NoTor:6} %
$K_0 (\Cs{2}{A}{\af})$ has no torsion,
and $K_1 (\Cs{2}{A}{\af}) \cong \Z.$
\end{enumerate}
\end{prp}

\begin{proof}
The algebra $A$ is the same as in Proposition~\ref{P:Torsion},
so~(\ref{P:NoTor:1}) and the first statement in~(\ref{P:NoTor:5})
follow from the corresponding parts of Proposition~\ref{P:Torsion}.
The rest of~(\ref{P:NoTor:5}),
the fact that $\af$ generates an action with the \tRp,
and the first part of~(\ref{P:NoTor:4}),
all have the same proofs as
the corresponding parts of Proposition~\ref{P:Torsion}.

The computation of the crossed products is now different.
Set
\[
E = \{ (x_0, x_1, \ldots, x_{2_m}) \in S^{2 m} \colon x_0 = 0 \}
\]
and
\[
L = \{ (x_0, x_1, \ldots, x_{2_m}) \in S^{2 m} \colon x_0 \geq 0 \}.
\]
Note that $E \cong S^{2 m - 1}.$
There is an injective homomorphism $\mu \colon B \to C (L, M_2)$
which sends $f \in C ( S^{2 m})$ to the function
$\mu (f) (x) = \diag ( f(x), f (h (x)))$ and the standard implementing
unitary in the crossed product to the constant function with value
\[
\left( \begin{array}{cc} 0 & 1 \\
1 & 0 \end{array} \right).
\]
Its image is
\[
\{ f \in C ( L, M_2) \colon {\mbox{$f (x)_{1, 1} = f (x)_{2, 2}$ and
$f (x)_{2, 1} = f (x)_{1, 2}$ for $x \in E$}} \}.
\]
Conjugating by the unitary
\[
\frac{1}{\sqrt{2}}
 \left( \begin{array}{cc} 1 & 1 \\ -1 & 1 \end{array} \right),
\]
we identify $B$ with
\[
\{ f \in C (L, M_2) \colon
   {\mbox{$f (x)$ is diagonal for $x \in E$}} \}.
\]

We compute $K_* (B)$ via a trick.
Set
\[
e = \left( \begin{array}{cc} 1 & 0 \\
0 & 0 \end{array} \right)
\andeqn
C = \{ f \in C ( L, M_2) \colon {\mbox{$f (x) \in \C e$ for
$x \in E$}} \}.
\]
The map $C (L) \to C,$ sending $f \in C (L)$ to the function
$x \to f (x) e,$ is an isomorphism on K-theory, essentially
because $K \otimes C \cong K \otimes C (L).$
Thus, $K_0 (C) \cong \Z,$ generated by the class of the constant
function ${\overline{e}}$ with value $e,$ so
$K_0 (C) \to K_0 (B)$ is injective.
There is a short exact sequence
\[
0 \longrightarrow C \longrightarrow B
 \longrightarrow C (S^{2 m - 1}) \longrightarrow 0,
\]
and the associated six term exact sequence in K-theory is
\[
\begin{CD}
\Z @>>> K_0 (B) @>>> \Z   \\
@AAA & &  @VVV            \\
\Z @<<< K_1 (B) @<<< \, 0.
\end{CD}
\]
We have just seen that the map $\Z \to K_0 (B)$ is injective, and
$K_0 (C (S^{2 m - 1})) \cong \Z,$ generated by the class of the
identity.
So $K_0 (B) \cong \Z^2,$ generated by $[ {\overline{e}} ]$ and
$[ 1 - {\overline{e}} ],$ and $K_1 (B) \cong \Z.$

As in the proof of Proposition~\ref{P:Torsion}, we have
$B_n = C^* ( \Z_2, A_n, \af_n ) \cong M_{s (n)} \otimes B.$
Moreover, by the same reasoning as there,
the map $B_{n - 1} \to B_n \cong M_{2 r (n) + 1} \otimes B_{n - 1}$
is the direct sum of the identity map and a map which factors
through a finite dimensional \ca.
Therefore $K_1 (B_{n - 1}) \to K_1 (B_n)$ is an isomorphism
for all $n,$
whence $K_1 (C^* (\Z_2, A, \af)) \cong \Z.$
Furthermore, $K_0 (C^* (\Z_2, A, \af))$ has no torsion because
it is a direct limit of groups isomorphic to $\Z^2.$
We now have~(\ref{P:NoTor:6}).
We also know that $C^* (\Z_2, A, \af)$ is not~AF,
so Theorem~2.2 of~\cite{PhtRp1a} implies that $\af$ does not have
the \sRp.
This finishes the proof of~(\ref{P:NoTor:3}).

We can also prove the second part of~(\ref{P:NoTor:4}).
Suppose $\af$ is \sar.
Then ${\widehat{\af}}$ has the \tRp\  %
by Theorem~\ref{RPDualToAppRep}(\ref{RPDualToAppRep:2}).
Using dual action, $K_1 (C^* (\Z_2, A, \af)) \cong \Z$ becomes a
module over the group ring $\Z [ \Z_2 ],$
with the nontrivial group element acting via the
nontrivial automorphism.
Theorem~3.3 of~\cite{Iz2} implies that
$K_1 (C^* (\Z_2, A, \af))$ is cohomologically trivial
in the sense of Definition~3.1(2) of~\cite{Iz2}.
But there are only two ways to make $\Z$ a $\Z [ \Z_2 ]$-module.
If the nontrivial group element acts trivially,
then, in the notation of~\cite{Iz2}, we have
${\mathrm{Coker}} \big( {\overline{N}} \big) \cong \Z_2,$
while if nontrivial group element acts by multiplication by $-1,$
then ${\mathrm{Ker}} \big( {\overline{N}} \big) \cong \Z_2.$
In either case,
$K_1 (C^* (\Z_2, A, \af))$ is not cohomologically trivial,
and this contradiction shows that $\af$ is not \sar.

It remains to prove~(\ref{P:NoTor:2}).
It follows from Theorem~2.6 of~\cite{PhtRp1a}
that $C^* (\Z_2, A, \af)$ has tracial rank zero.
We have seen that it is a direct limit of type~I C*-algebras, so
that it is nuclear and satisfies the Universal Coefficient Theorem.
It is simple by Corollary~1.6 of~\cite{PhtRp1a}.
Therefore Lin's classification theorem
(Theorem~5.2 of~\cite{Ln15}) applies.
In the form given in Proposition~3.7 of~\cite{PhtRp2},
it implies that $C^* (\Z_2, A, \af)$
is a simple AH~algebra with real rank zero
and no dimension growth.
Proposition~5.7 of~\cite{ELP} implies that
$C^* (\Z_2, A, \af)$ has a unique tracial state.
This completes the proof.
\end{proof}

Unlike in the situation of Example~\ref{E:Torsion},
the classification theorem of~\cite{EGL} does not apply
to the direct limit,
and we do not know of a classification proof which does not
depend on the \tRp.

\begin{rmk}\label{R:Inf}
In Examples~\ref{E:Torsion} and~\ref{E:NoTor}, the requirement
$\lim_{n \to \I} r (n) = \I$ is merely for convenience in proving the
\tRp.
All features, including the \tRp,
hold provided $r (n) \geq 1$ for all $n.$
In Example~\ref{E:NoTor} (but not in Example~\ref{E:Torsion}),
this increases the number of possibilities
for $K_0 (C^* ( \Z_2, A, \af)).$
\end{rmk}

\begin{rmk}\label{R:Even}
In Example~\ref{E:NoTor}, we do not need to require that all the
matrix sizes be odd.
For example, we could use in place of $\ps_{r (n), \, x (n)}$ the map
which sends $f \in C (S^{2 m})$ to
\[
x \mapsto \diag ( f(x), \, f(y (n)), \, f (z (n)), \, f (h (z (n))))
\]
with $y (n), z (n) \in S^{2 m}$ chosen such that $h (y (n)) = y (n)$ and
$h (z (n)) \neq z (n).$
In place of $w_{r (n)}$ one uses
\[
\left( \begin{array}{cccc} 1 & 0 & 0 & 0 \\
0 & 1 & 0 & 0 \\
0 & 0 & 0 & 1 \\
0 & 0 & 1 & 0 \end{array} \right).
\]
One gets an action of $\Z_2$ on the AF~algebra in~10.11.3
of~\cite{Bl} which induces the map there on K-theory.

This particular variation does not work in Example~\ref{E:Torsion},
because there are no fixed points.
\end{rmk}

\begin{rmk}\label{R:NoRokhlin}
Let $A$ be as in Remark~\ref{R:Even}.
As in~10.11.3 of~\cite{Bl}, there is no $\af \in \Aut (A)$
satisfying the following:
\begin{enumerate}
\item\label{R:NoRokhlin:1}
$\af^2 = \id_A.$
\item\label{R:NoRokhlin:2}
$\af_* ( \et, k) = ( \et, \, - k)$
on $K_0 (A) = \Z \big[ \frac{1}{2} \big] \oplus \Z.$
\item\label{R:NoRokhlin:3}
$C^* ( \Z_2, A, \af)$ is~AF.
\end{enumerate}
By Theorem~2.2 of~\cite{PhtRp1a},
there is therefore no $\af \in \Aut (A)$
satisfying~(\ref{R:NoRokhlin:1}) and~(\ref{R:NoRokhlin:2}) and
such that $\af$ generates an action of $\Z_2$
with the \sRp.
This is a less obvious nonexistence result than for, say, actions
of $\Z_2$ on the $3^{\I}$~UHF algebra.
However, it can also be obtained from Theorem~3.3 of~\cite{Iz2}.
\end{rmk}

\section{Questions}\label{Sec:Qst}

\indent
The \tRp, as given in Definition~1.2 of~\cite{PhtRp1a},
is certainly not useful for \ca s with few \pj s,
may well not be useful for \ca s with many \pj s but bad
comparison theory,
and may not be useful for nonsimple \ca s.
The same is true of tracial approximate representability
and tracial approximate innerness (Definition~5.1 of~\cite{PhtRp1a}).
In this section,
we give a brief discussion of what might be done for the \tRp.
Similar considerations apply to the other properties,
but we have less information on which to base any speculation.

If the \ca\  $A$ is purely infinite,
then the condition in Lemma~\ref{TRPCond}
for an action $\af \colon G \to \Aut (A)$ to have the
\tRp\  is vacuous,
since $e_g = 0$ for all $g \in G$
still permits Condition~(\ref{TRPCond:3}) to hold.
Definition~1.2 of~\cite{PhtRp1a} contains the following extra
condition:
\begin{enumerate}
\setcounter{enumi}{\value{TmpEnumi}}
\item\label{TRPCond:4} %
With $e = \sum_{g \in G} e_g$ (as in Condition~(\ref{TRPCond:3})),
we have $\| e x e \| > 1 - \ep.$
\setcounter{TmpEnumi}{\value{enumi}}
\end{enumerate}
(This condition was inspired by Definition~2.1 of~\cite{LnTAF}.)
The definition we originally gave in~\cite{PhW}
required, instead of Condition~(\ref{TRPCond:4}),
that for every $N \in \N$
one be able to choose the \pj s $e_g$ so that, in addition,
\begin{enumerate}
\setcounter{enumi}{\value{TmpEnumi}}
\item\label{TRPCond:5} %
With $e = \sum_{g \in G} e_g$ as before,
there are $N$ \mops\  $f_1, f_2, \dots, f_N \leq e_g,$
each of which is \mvnt\  to $1 - e.$
\end{enumerate}
Under either extra condition,
one can force $e_g \neq 0$ for all $g \in G.$
The proof of Lemma~1.14 of~\cite{PhtRp1a}
then shows that $\af_g$ is outer for all $g \in G \setminus \{ 1 \}.$

According to Theorem~1 of~\cite{Nk},
if $A$ is a unital Kirchberg algebra
(purely infinite simple separable and nuclear),
and if $\af \in \Aut (A)$ has the property that $\af^n$ is outer
for every $\af \neq 0,$
then $\af$ has the \sRp.
For finite group actions, there are obvious K-theoretic
obstructions to the \sRp.
For example, if $A$ is any unital purely infinite simple \ca\  with
$K_0 (A) \cong \Z$ and generated by $[1],$
such as ${\mathcal{O}}_{\infty},$
then every automorphism is trivial on $K_0.$
If $\af \colon G \to \Aut (A)$ had the \sRp,
then $[1] \in K_0 (A)$ would then be divisible by $\card (G),$
forcing $G$ to be trivial.
For classifiable algebras there are less obvious K-theoretic
obstructions, such as Theorem~3.5 of~\cite{Iz2}.
However, the following conjecture seems plausible.

\begin{cnj}\label{PI}
Let $A$ be a unital Kirchberg algebra,
and let $\af \colon G \to \Aut (A)$ be an action of a finite group $G$
on $A$
such that $\af_g$ is outer for all $g \in G \setminus \{ 1 \}.$
Then $\af$ has the \tRp,
as given in Definition~1.2 of~\cite{PhtRp1a}.
\end{cnj}

Since the proof in~\cite{Nk} uses nuclearity in an
apparently essential way,
we do not know what to expect in the nonnuclear case.
For algebras which are infinite but not purely infinite,
the situation is even less clear.

We now turn to algebras with few \pj s.
The following discussion applies only to simple unital \ca s
with ``enough traces'';
if, for example, the algebra is stably finite but the
comparison theory of \pj s is bad,
the smallness conditions on the leftovers must be changed.

The criterion for the \tRp\  in Theorem~\ref{ARPFromPosElts2}
does not mention \pj s,
and thus makes sense even if there are no nontrivial \pj s.
We suppose, for example, that it probably holds for the
flip $a \otimes b \mapsto b \otimes a$ on the tensor
product of two copies of the Jiang-Su algebra $Z$~\cite{JS}.
(As shown in~\cite{OP3}, for many simple \ca s which have many \pj s,
the flip has the \tRp.)
This suggests taking the criterion in Theorem~\ref{ARPFromPosElts2}
as the definition in the absence of enough \pj s.
Since Theorem~\ref{ARPFromPosElts2} depends on tracial rank zero,
this suggestion may need to be limited to algebras which are expected
to be classifiable.

There are two tests.
First, one might hope for an analog of Theorem~2.6 of~\cite{PhtRp1a}
and its generalization in~\cite{OP3}.
To include algebras with few \pj s,
and lacking an abstract criterion for classifiability for such algebras,
we propose the following problem.

\begin{pbm}\label{RSHA}
Let $A$ be an \idsuca\  which is isomorphic to a
direct limit, with no dimension growth,
of recursive subhomogeneous algebras,
as in Theorem~3.6 of~\cite{Ph1}.
Let $\af \colon G \to \Aut (A)$
be an action of a finite group $G$ on $A,$
and suppose that $\af$ satisfies the condition of
Theorem~\ref{ARPFromPosElts2}.
Does it follow that $C^* (G, A, \af)$
is again isomorphic to a direct limit, with no dimension growth,
of recursive subhomogeneous algebras?
\end{pbm}

Note that the class contains the Jiang-Su algebra $Z,$
as well as other simple unital \ca s with no nontrivial \pj s.

Of course, if the answer is no,
this might mean that the class of \ca s is wrong,
rather than indicating any any problem
with the proposed definition of the \tRp.

For the second test, we begin with a conjecture which seems
very plausible.

\begin{cnj}\label{RankOfCP}
Let $A$ be an \idsuca, let $G$ be a finite group,
and let $\af \colon G \to \Aut (A)$ be an action with the \tRp.
Suppose that $A$ has real rank zero, stable rank one,
and that the order on \pj s over $A$ is determined by traces
(Definition~2.4 of~\cite{PhtRp1a};
a stable version of Blackadar's Second Fundamental Comparability Question
as in 1.3.1 in~\cite{Bl3}).
Then $C^* (G, A, \af)$ also has these three properties.
\end{cnj}

This is true for actions of $\Z$ on stably finite \suca s,
by Theorems~3.5, 4.5, and~5.3 of~\cite{OP1}.
The proofs for actions of finite groups should be similar but simpler.
The expected truth of Conjecture~\ref{RankOfCP}
motivates the following problem.
We drop real rank zero,
since we are interested in algebras with few \pj s.

\begin{pbm}\label{RankPbm}
Let $A$ be an \idsuca, let $G$ be a finite group,
and let $\af \colon G \to \Aut (A)$ be an action
which satisfies the condition of Theorem~\ref{ARPFromPosElts2}.
Suppose that $A$ has stable rank one
and that the order on \pj s over $A$ is determined by traces.
Does $C^* (G, A, \af)$ also have these two properties?
\end{pbm}

We now turn to nonsimple \ca s.
Ultimately,
one hopes for a general theory encompassing actions on simple \ca s,
on commutative \ca s,
and on everything in between.
It is not clear what the correct definitions are.
The commutative case with many \pj s is not very informative,
because of the following result.

\begin{prp}\label{Cantor}
Let the finite group $G$ act on the Cantor set $X,$
and let $\af \colon G \to \Aut (C (X))$
be the action $\af_g (f) (x) = f (g^{-1} x)$
for $x \in X,$ $g \in G,$ and $f \in C (X).$
Then \tfae:
\begin{enumerate}
\item\label{Cantor:1} %
$\af$ has the \sRp.
\item\label{Cantor:2} %
$\af$ satisfies the condition of the definition of the \tRp\  %
(Definition~1.2 of~\cite{PhtRp1a}; see the beginning of this section).
\item\label{Cantor:3} %
$\af$ satisfies the condition of Lemma~\ref{TRPCond}.
\item\label{Cantor:4} %
The action of $G$ on $X$ is free.
\end{enumerate}
\end{prp}

\begin{proof}
It is clear that (\ref{Cantor:1}) implies~(\ref{Cantor:2})
and (\ref{Cantor:2}) implies~(\ref{Cantor:3}).

Assume~(\ref{Cantor:3}).
Let $x \in X.$
Choose a nonempty compact open set $K \S X$
such that $K \cap G x = \varnothing.$
Apply~(\ref{Cantor:3}) with $F = \varnothing,$
with $\ep = \frac{1}{2},$
and with $x = \ch_K,$
obtaining disjoint compact open sets $L_g$ for $g \in G$ such that
$\| \ch_{g L_h} - \ch_{L_{g h}} \| < \frac{1}{2}$ for all $g, h \in G,$
and such that, with $L = \bigcup_{g \in G} L_g,$
we have $\ch_{X \setminus L} \precsim \ch_K.$
It follows that $g L_h = L_{g h}$ for all $g, h \in G,$
and that $X \setminus L \subset K.$
Since $x \in L,$ the map $g \mapsto g x$ is injective.
Since this is true for all $x \in X,$ we have proved~(\ref{Cantor:4}).

Now assume~(\ref{Cantor:4}).
For each $x \in X,$ use continuity and freeness to choose
a compact open set $K_x \subset X$ which contains $x$ and such that
the sets $g K_x,$ for $g \in G,$ are disjoint.
Choose $x_1, x_2, \ldots, x_n \in X$ such that
$K_{x_1}, K_{x_2}, \ldots, K_{x_n}$ cover $X.$
Inductively set $N_1 = K_{x_1}$ and
$N_{k + 1} = N_k \cup [K_{x_{k + 1}} \cap (X \SM G N_k)].$
The sets $g N_1,$ for $g \in G,$ are disjoint.
Moreover, if the sets $g N_k,$ for $g \in G,$ are disjoint,
then for distinct $g, h \in G$ we have
\begin{align*}
g N_{k + 1} \cap h N_{k + 1}
& \subset [g N_k \cap h N_k]
        \cup [g N_k \cap h (X \SM G N_k)]   \\
& \hspace*{3em}        \cup [g (X \SM G N_k) \cap h N_k]
        \cup [g K_{x_{k + 1}} \cap h K_{x_{k + 1}}].
\end{align*}
All four terms are clearly empty.
By induction, therefore, the sets $g N_n,$ for $g \in G,$ are disjoint.
A related argument shows that $G N_k$ contains
$G K_{x_1}, G K_{x_2}, \ldots, G K_{x_n},$
whence $G N_n = X.$

Now the \pj s $e_g = \ch_{g N_n}$ are central,
satisfy $\sum_{g \in G} e_g = 1,$
and satisfy $\af_g (e_h) = e_{g h}$ for all $g, h \in G.$
Accordingly, $\af$ has the \sRp.
\end{proof}

We suspect, however,
that the condition of Lemma~\ref{TRPCond},
or even the \tRp\  as defined here,
will not be adequate to prove general theorems
about actions on nonsimple \ca s.

\end{document}